\documentclass[a4paper,12pt]{amsart}

\usepackage{amsmath,amsthm,amssymb}
\usepackage[all]{xy}

\swapnumbers
\theoremstyle{plain}
  \newtheorem{thm}{Theorem}[section]
  \newtheorem*{thm*}{Theorem}
  \newtheorem{lem}[thm]{Lemma}
  \newtheorem*{lem*}{Lemma}
  \newtheorem{cor}[thm]{Corollary}
  \newtheorem*{cor*}{Corollary}
  \newtheorem{prop}[thm]{Proposition}
  \newtheorem*{prop*}{Proposition}
\theoremstyle{definition}

\theoremstyle{remark}
  \newtheorem{rem}[thm]{Remark}
  \newtheorem*{rem*}{Remark}

\newcommand{\nn}[1]{(\ref{#1})}
\newcommand{\cP}{{\Cal P}}
\newcommand{\ol}[1]{\overline{#1}}\newcommand{\ul}[1]{\underline{#1}}

\def\Cal{\mathcal}

\newcommand{\cq}{{\Cal Q}}
\newcommand{\ce}{{\Cal E}}
\newcommand{\bg}{\mbox{\boldmath{$ g$}}}
\newcommand{\nd}{\nabla}
\newcommand{\Rho}{P}%{{\mbox{\sf P}}}
\newcommand{\Ric}{\operatorname{Ric}}
\newcommand{\J}{J}%{{\mbox{\sf J}}}
\newcommand{\V}{P}%{{\mbox{\sf P}}}    
\newcommand{\wh}{\widehat}
\newcommand{\wt}{\widetilde}

\def\al{\alpha}

\def\de{\delta}

\def\ka{\kappa}

\def\rh{\rho}
\def\si{\sigma}
\def\ta{\tau}
\def\ph{\varphi}

\def\om{\omega}

\def\De{\Delta}

\def\Ph{\Phi}
\def\Ps{\Psi}
\def\Om{\Omega}
\def\Up{\Upsilon}
\def\na{\nabla}

\def\D{\mathbb{D}}

\def\W{\mathbb{W}}
\def\X{\mathbb{X}}
\def\Y{\mathbb{Y}}
\def\Z{\mathbb{Z}}

\def\cE{\mathcal{E}}

\def\InR{\in \mathbb{R}}

\def\ekv{\Longleftrightarrow}
\def\form#1{\mathbf{#1}}
\def\dform#1{\dot{\mathbf{#1}}}
\def\ddform#1{\ddot{\mathbf{#1}}}
\def\dddform#1{\dddot{\mathbf{#1}}}

\def\idx#1{{\em #1\/}}
\renewcommand{\vec}[1]{\mathbf{#1}}

\newcommand{\lpl}{
  \mbox{$
  \begin{picture}(12.7,8)(-.5,-1)
  \put(2,0.2){$+$}
  \put(6.2,2.8){\oval(8,8)[l]}
  \end{picture}$}}

\def\sideremark#1{\ifvmode\leavevmode\fi\vadjust{\vbox to0pt{\vss% the remark
 \hbox to 0pt{\hskip\hsize\hskip1em%                          will appear only
 \vbox{\hsize3cm\tiny\raggedright\pretolerance10000%          on the side
 \noindent #1\hfill}\hss}\vbox to8pt{\vfil}\vss}}}%
\begin{document}

\author{A. Rod Gover and Josef \v Silhan} \title{The conformal Killing
  equation on forms -- prolongations and applications} \date{}

\begin{abstract}
  We construct a conformally invariant vector bundle connection such
  that its equation of parallel transport is a first order system that
  gives a prolongation of the conformal Killing equation on
  differential forms. Parallel sections of this connection are related
  bijectively to solutions of the conformal Killing equation.  We
  construct other conformally invariant connections, also giving
  prolongations of the conformal Killing equation, that bijectively
  relate solutions of the conformal Killing equation on $k$-forms to a
  twisting of the conformal Killing equation on $(k - \ell)$-forms
  for various integers $\ell$. These tools are used to develop a helicity
  raising and lowering construction in the general setting and on 
conformally Einstein manifolds.
\end{abstract}

\maketitle

\pagestyle{myheadings}
\markboth{Gover \& \v Silhan}{Conformal Killing equations}

\section{Introduction}

On a (pseudo-)Riemannian $n$-manifold a tangent vector field $v$ is an
infinitesimal conformal automorphisms if the Lie derivative of the
metric ${\Cal L}_v g$ is proportional to the metric $g$.  This is the
so-called conformal Killing equation and, denoting by $\si$ the 1-form
$g(v,~)$, the equation may be re-expressed as the equation requiring
the trace-free symmetric part of $\na \si$ to be zero.  
Here we will use the term conformal Killing equation to mean the first
 order overdetermined partial differential
 equation which generalises this to (conformally weighted)
 differential forms of rank $k$ where $1\leq k \leq n-1$: a $k$-form $\si$ 
is a conformal Killing form if, with respect to the $O(g)$-decomposition of $T^*M\otimes \Lambda^k T^*M$, the Cartan part of $\na \si$ is zero. Equivalently for any tangent vector field $u$ we have 
\begin{equation}\label{ifv}
\na_u \si= \varepsilon(u)\tau + \iota(u) \rho
\end{equation}
where, on the right-hand side $\tau$ is a $k-1$-form, $\rho$ is a
$k+1$-form, and $\varepsilon(u)$ and $\iota(u)$ indicate,
respectively, the exterior multiplication and (its formal adjoint) the
interior multiplication of $g(u,~)$. An important property of the
conformal Killing equation \nn{ifv} is that it is conformally
invariant (where we require the $k$-form $\si$ to have conformal
weight $k+1$).  

In the simplest terms the main aims of this paper are:
1.\ to derive a conformally invariant connection which is
``equivalent'' to the conformal Killing equation  in the sense that its
parallel sections are naturally in one-one correspondence with
solutions of the of the conformal Killing equation \nn{ifv}, 2.\ to derive a
conformally invariant connection $\tilde{\na}$ that in a similar way
relates solutions of the of the conformal Killing equation on
$k$-forms to solutions of the solutions of the of the conformal
Killing equation on $(k-\ell)$-forms (for suitable positive and negative integers $\ell$) twisted by the connection
$\tilde{\na}$, and 3.\ to apply these ideas to a programme of helicity
raising and lowering (in the sense of \cite{nt}, for example in the
``helicity raising'' direction, given certain assumptions of the
conformal curvature, pairs of conformal Killing forms may be ``cupped''
together to yield a conformal Killing form of higher rank). For the third part
here no attempt has been made to be complete. Rather our philosophy
has been to establish some basic results in this direction, which
follow easily from the machinery established in parts 1.\ and 2.,\ and
through this indicate the broad idea and explore some applications.

Conformal Killing 2-forms were introduced by Tachibana in \cite{Tach}
and the generalisation to higher valence followed shortly after
\cite{Kashiwada}. Coclosed conformal Killing forms are Killing forms (or
sometime called Killing-Yano forms). The latter satisfy the equation
which generalises the Killing equation on vector fields, that is
\nn{ifv} with $\tau$ identically 0. This equation has been studied
intensively in the Physics literature in connection with its role
generating quadratic first integrals of the geodesic equation. Aside
from this connection, and a role in the higher symmetries of other
equations \cite{BennC}, the broader geometric meaning of higher rank
conformal Killing forms is still somewhat mysterious. The linear
operator giving the conformal Killing equation is a Stein-Weiss
gradient and elliptic in the Riemannian setting \cite{B-StW}. The
issue of global existence of conformal Killing forms in the compact
Riemannian setting has been pursued recently by Semmelmann and others,
see \cite{Sem-MathZ,SemmG2} for an indication of results and further references. Our
treatment here is primarily in arbitrary signature and concerned with
local issues. In particular we seek to draw out the additional
information arising from the conformal invariance of the equation.
This should have important consequences for the general theory
including global existence.

It is often the case that an overdetermined linear partial
differential equation is equivalent to a first order prolonged system
that may be interpreted as a vector bundle connection and its equation
of parallel transport.  More generally a semilinear partial
differential equation is said to be of {\em finite type} if there is a
suitably equivalent finite dimensional prolonged system that is
``closed'' in the sense that all first partial derivatives of the
dependent variables are determined by algebraic formulae in terms of
these same variables. There is a criteria due to Spencer
\cite{spencer} to determine when a semilinear equation has this
property. However for any given equation one generally wants
considerably 
more information. For the case of the conformal Killing
equation Semmelmann explicitly constructs a prolongation and
connection along these lines \cite{Sem-MathZ} and so establishes sharp
bounds on the dimension of spaces of conformal Killing forms. This
idea was generalised in \cite{BCEG} where ideas from Kostant's Hodge
theory are used to give a uniform algorithm for explicitly computing
such prolongations for a wide class of geometric semilinear
overdetermined partial differential equations. This class includes the
conformal Killing equation as one of the simplest cases. However
neither of these treatments addresses the conformal invariance of
conformal Killing equation. For the case of conformal Killing
equations on vector fields an equivalent conformally invariant
connection was given in \cite{powerslap}. (See also \cite{Andi-infint}
which generalises this by giving an invariant connection corresponding
to the equation of infinitesimal automorphisms for all parabolic
geometries.)  Ab initio, given a conformally invariant equation one
does not know whether there is a conformally invariant prolonged
system along these lines. We show that for the conformal Killing
equation there is, see Theorem \ref{main}. In fact this theorem gives
much more, it gives a connection equivalent (in the sense of the
theorem) to the conformal Killing equation and this conformally
invariant connection is described explicitly in terms of the normal
tractor connection of \cite{BEGo,Cap+Gover}. The power of this is that
the normal tractor connection (here on a bundle of rank
$(n+1)(n+2)/2$) is an exterior tensor product of the (normal)
connection on standard tractors, and the latter is a simple well-understood 
connection on a bundle of relatively low rank (viz.\ $n+2$) and
which respects a bundle metric. By describing things in terms of {\em
form-tractors} in this way we also capture succinctly what conformal
invariance means for the components of the prolongation. (The form
form-tractor bundles and their normal connection are treated
explicitly in \cite{Branson+Gover}.)

An application of Theorem \ref{main} (that we sketch here but shall not pursue
 explicitly in this work) is that it provides a way to construct
 natural conformal invariants (conformal invariants that may be
 expressed by universal formulae involving complete or partial
 contractions of the Riemannian curvature and its covariant
 derivatives) which locally obstruct the existence of conformal
 Killing forms. Since such forms correspond to a holonomy reduction of
 the connection obtained in the Theorem, it is clear that its
 curvature will in general obstruct such forms. By analogy with the
 treatment of obstructions to conformally Einstein metrics in Section 3.3 of
 \cite{GoNur}, an appropriately defined ``determinant'' of this curvature
 must vanish in order for the conformal Killing equation to have a
 solution. From this it is straightforward to extract the required natural
 invariants.

We may view a connection as a twisting (or coupling) of the exterior
derivative on functions. The latter is of course conformally
invariant.  Thus Theorem \ref{main} relates solutions of the conformal
Killing equation to a conformal twisting of the exterior derivative. A
generalisation of this idea (which to our knowledge has not been
explored previously) is, given two suitable distinct conformally
invariant equations $A=0$ and $B=0$, to consider obtaining a
conformally connection $\tilde{\nd}$ so that solutions of the equation
$A$ are bijectively related (by a prolongation) to solutions of the
twisting by $\tilde\nd$ of the equation $B$ i.e.\
$B^{\tilde{\nd}}=0$. (Of course there are variants of this where we
replace conformal invariance by any other notion of invariance. We
should also point out that this idea is related to, and in special
cases agrees with, the ``translation principle'' touched on in section
8 of \cite{CaldD}.) In Section \ref{s.helicity} we obtain results
exactly of this type, with the conformal Killing equation on forms of
different ranks playing the roles of both $A$ and $B$, see Theorem
\ref{coupKillthm}, and also Proposition \ref{coupKill} which deals
case that ``$B$'' is the conformal Killing equation on tangent
vectors. The remarkable feature of these results is the very simple
form of the twisting connection -- it differs from the normal tractor
connection by a simple curvature action as in expressions
\nn{mod_connection} and \nn{dcoup}.  These results and their
simplicity are exploited in Section \ref{almosthell} where we describe
explicit conformally invariant helicity raising and lowering formulae
and their obstructions to being a solution of the conformal Killing
equation.  See in particular: Theorem \ref{t.einstein} which uses 
(almost) Einstein metrics and conformal Killing fields to generate conformal Killing fields; 
 Theorem \ref{t.ckf} where conformal Killing
forms are used to generate other conformal Killing forms; and Theorem
\ref{t.ckf} where they are used to generate conformal Killing {\em
tensors}, i.e. symmetric trace-free tensors $S$ such that the
symmetric trace-free part of $\nd S$ vanishes.
This idea of combining solutions to yield solutions of other equations is
along the lines of helicity raising and lowering by Penrose's twistors
in dimension 4. See also \cite{CaldD} where a helicity raising and
lowering machinery is developed, in very general terms, via
tractor-twisted differential forms and the twisted exterior
derivative.

To construct the required prolongations we develop and employ an
efficient calculus that enables us to efficiently deal with
differential forms and form-tractors and some related bundles of
arbitrary rank. Some of the ideas for this originate in
\cite{Branson+Gover} but many new techniques and tools have been
developed in \cite{josth} and the reader can find greater on this (and
many aspects of this article) in that source. Very crudely the idea is
that as a first step in constructing the new connections we may take
the normal tractor connection (or its coupling with the Lev-Civita
connection) to be a ``first approximation'' to the required new
connection. By elementary representation theory it in fact must agree
with the required connections in the conformally flat setting.  Then,
employing the form calculus mentioned, we compute explicitly the
tractor ``contorsion'' needed to adjust the normal connection.
Eastwood's curved translation principle (see e.g\ \cite{MikESrni}) is
a technique for generating conformally invariant equations from other
such equations via differential splitting operators between tractor
bundles and (weighted) tensor-spinor bundles. The constructions and
ideas in sections \ref{s.helicity} and \ref{almosthell} are partly
inspired by this technique and involve the refinement where one seeks to
``translate solutions'' of equations rather just the equations themselves,
this necessarily draws on solutions of other equations and their
equivalence to parallel (or suitably almost parallel) sections of
tractor bundles, thus a form of helicity raising-lowering.  

We would
like to thank Mike Eastwood and Andi \v Cap for discussions concerning
their views on the prolongation treatment of conformal Killing vectors
and related issues.

\section{Conformal geometry, tractor calculus and conformal Killing equation}

\subsection{Conformal geometry and tractor calculus}

We summarise here some relevant notation and background
for conformal structures.  Further details may be found in
\cite{CapGoamb,GoPetLap}.  Let $M$ be a smooth manifold of dimension
$n\geq 3$. Recall that a {\em conformal structure\/} of signature
$(p,q)$ on $M$ is a smooth ray subbundle $\cq\subset S^2T^*M$ whose
fibre over $x$ consists of conformally related signature-$(p,q)$
metrics at the point $x$. Sections of $\cq$ are metrics $g$ on $M$. So
we may equivalently view the conformal structure as the equivalence
class $[g]$ of these conformally related metrics.  The principal
bundle $\pi:\cq\to M$ has structure group $\Bbb R_+$, and so each
representation ${\Bbb R}_+ \ni x\mapsto x^{-w/2}\in {\rm End}(\Bbb R)$
induces a natural line bundle on $ (M,[g])$ that we term the conformal
density bundle $E[w]$. We shall write $ \ce[w]$ for the space of
sections of this bundle.  We write $\ce^a$ for the space of sections
of the tangent bundle $TM$ and $\ce_a$ for the space of sections of
$T^*M$. The indices here are abstract in the sense of \cite{ot} and we
follow the usual conventions from that source. So for example
$\ce_{ab}$ is the space of sections of $\otimes^2T^*M$.  Here and
throughout, sections, tensors, and functions are always smooth.  When
no confusion is likely to arise, we will use the same notation for a
bundle and its section space.

We write $\bg$ for the {\em conformal metric}, that is the
tautological section of $S^2T^*M\otimes E[2]$ determined by the
conformal structure. This will be used to identify $TM$ with
$T^*M[2]$.  For many calculations we will use abstract indices in an
obvious way.  Given a choice of metric $ g$ from the conformal class,
we write $ \nabla$ for the corresponding Levi-Civita connection. With
these conventions the Laplacian $ \Delta$ is given by
$\Delta=\bg^{ab}\nd_a\nd_b= \nd^b\nd_b\,$. Here we are raising indices
and contracting using the (inverse) conformal metric. Indices will be
raised and lowered in this way without further comment.  Note $E[w]$
is trivialised by a choice of metric $g$ from the conformal class, and
we write $\nd$ for the connection corresponding to this
trivialisation.  It follows immediately that (the coupled) $ \nd_a$
preserves the conformal metric.

Since the Levi-Civita connection is torsion-free, its curvature
$R_{ab}{}^c{}_d$ (the Riemannian
curvature) is given by $ [\nd_a,\nd_b]v^c=R_{ab}{}^c{}_dv^d $
($[\cdot,\cdot]$ indicates the commutator bracket).  The Riemannian
curvature can be decomposed into the totally trace-free Weyl curvature
$C_{abcd}$ and a remaining part described by the symmetric {\em
Schouten tensor} $\Rho_{ab}$, according to 
\begin{equation}\label{csplit}
R_{abcd}=C_{abcd}+2\bg_{c[a}\Rho_{b]d}+2\bg_{d[b}\Rho_{a]c}, 
\end{equation}
 where
$[\cdots]$ indicates antisymmetrisation over the enclosed indices.
The Schouten tensor is a trace modification of the Ricci tensor
$\Ric_{ab}=R_{ca}{}^c{}_b$ 
and vice versa: $\Ric_{ab}=(n-2)\Rho_{ab}+\J\bg_{ab}$,
where we write $ \J$ for the trace $ \V_a{}^{a}$ of $ \V$.  The {\em
Cotton tensor} is defined by 
$$
A_{abc}:=2\nabla_{[b}\Rho_{c]a} .
$$
Via the Bianchi identity this is related to the divergence of the Weyl 
tensor as follows:
\begin{equation}\label{bi1} (n-3)A_{abc}=\nabla^d
C_{dabc} . 
\end{equation} 
Under a {\em conformal transformation} we replace a choice of metric 
$g$ by the metric $\hat{g}=e^{2\Up} g$, where $\Up$ 
 is a smooth
function. We recall that, in particular, the Weyl curvature is
conformally invariant $\widehat{C}_{abcd}=C_{abcd}$. (Note there that
as a type $(0,4)$-tensor-density the Weyl curvature is a
density-valued tensor of conformal weight 2.) On the other hand the
Schouten tensor transforms according to
\begin{equation}\label{Rhotrans}
\textstyle \widehat{\V}_{ab}=\V_{ab}-\nd_a \Up_b +\Up_a\Up_b
-\frac{1}{2} \Up^c\Up_c \bg_{ab} 
\end{equation}
where $\Up_a = \na_a \Up$. 

Explicit formulae for the corresponding transformation of
the Levi-Civita connection and its curvatures are given in e.g.\ 
\cite{BEGo,GoPetLap}. From these, one can easily compute the 
transformation for a general valence (i.e.\ rank) $s$ section 
$f_{bc \cdots d} \in \cE_{bc \cdots d}[w]$ using the Leibniz rule:
\begin{equation} \label{grad_trans_gen}
\begin{split}
  \hat{\na}_{\bar{a}} f_{bc \cdots d}
  =&  \na_{\bar{a}} f_{bc \cdots d} + (w-s) \Up_{\bar{a}} f_{bc \cdots d}
     -\Up_{b} f_{\bar{a}c \cdots d} \cdots
     -\Up_{d} f_{bc \cdots \bar{a}} \\
   & +\Up^p f_{p c \cdots d} \bg_{b\bar{a}} \cdots
     +\Up^p f_{bc \cdots p} \bg_{d\bar{a}}.
\end{split}
\end{equation}

We next define the standard tractor bundle over $(M,[g])$.
It is a vector bundle of rank $n+2$ defined, for each $g\in[g]$,
by  $[\ce^A]_g=\ce[1]\oplus\ce_a[1]\oplus\ce[-1]$. 
If $\wh g=e^{2\Up}g$, we identify  
 $(\alpha,\mu_a,\tau)\in[\ce^A]_g$ with
$(\wh\alpha,\wh\mu_a,\wh\tau)\in[\ce^A]_{\wh g}$
by the transformation
\begin{equation}\label{transf-tractor}
 \begin{pmatrix}
 \wh\alpha\\ \wh\mu_a\\ \wh\tau
 \end{pmatrix}=
 \begin{pmatrix}
 1 & 0& 0\\
 \Up_a&\delta_a{}^b&0\\
- \tfrac{1}{2}\Up_c\Up^c &-\Up^b& 1
 \end{pmatrix} 
 \begin{pmatrix}
 \alpha\\ \mu_b\\ \tau
 \end{pmatrix} .
\end{equation}
It is straightforward to verify that these identifications are
consistent upon changing to a third metric from the conformal class,
and so taking the quotient by this equivalence relation defines the
{\em standard tractor bundle} $\ce^A$ over the conformal manifold.
(Alternatively the standard tractor bundle may be constructed as a
canonical quotient of a certain 2-jet bundle or as an associated
bundle to the normal conformal Cartan bundle \cite{luminy}.) On a
conformal structure of signature $(p,q)$, the bundle $\ce^A$ admits an
invariant metric $ h_{AB}$ of signature $(p+1,q+1)$ and an invariant
connection, which we shall also denote by $\nabla_a$, preserving
$h_{AB}$. Up up to isomorphism this the unique {\em normal conformal
tractor connection} \cite{Cap+Gover} and it induces a normal
connection on $\bigotimes \ce^A$ that we will also denoted by $\na_a$
and term the (normal) tractor connection.  In a conformal scale $g$,
the metric $h_{AB}$ and $\na_a$ on $\ce^A$ are given by
\begin{equation}\label{basictrf}
 h_{AB}=\begin{pmatrix}
 0 & 0& 1\\
 0&\bg_{ab}&0\\
1 & 0 & 0
 \end{pmatrix}
\text{ and }
\nabla_a\begin{pmatrix}
 \alpha\\ \mu_b\\ \tau
 \end{pmatrix}
 =
\begin{pmatrix}
 \nabla_a \alpha-\mu_a \\
 \nabla_a \mu_b+ \bg_{ab} \tau +\Rho_{ab}\alpha \\
 \nabla_a \tau - \Rho_{ab}\mu^b  \end{pmatrix}. 
\end{equation}
It is readily verified that both of these are conformally well-defined,
i.e., independent of the choice of a metric $g\in [g]$.  Note that
$h_{AB}$ defines a section of $\ce_{AB}=\ce_A\otimes\ce_B$, where
$\ce_A$ is the dual bundle of $\ce^A$. Hence we may use $h_{AB}$ and
its inverse $h^{AB}$ to raise or lower indices of $\ce_A$, $\ce^A$ and
their tensor products.

In computations, it is often useful to introduce 
the `projectors' from $\ce^A$ to
the components $\ce[1]$, $\ce_a[1]$ and $\ce[-1]$ which are determined
by a choice of scale.
They are respectively denoted by $X_A\in\ce_A[1]$, 
$Z_{Aa}\in\ce_{Aa}[1]$ and $Y_A\in\ce_A[-1]$, where
 $\ce_{Aa}[w]=\ce_A\otimes\ce_a\otimes\ce[w]$, etc.
 Using the metrics $h_{AB}$ and $\bg_{ab}$ to raise indices,
we define $X^A, Z^{Aa}, Y^A$. Then we
immediately see that 
$$
Y_AX^A=1,\ \ Z_{Ab}Z^A{}_c=\bg_{bc} ,
$$
and that all other quadratic combinations that contract the tractor
index vanish. 
In \eqref{transf-tractor} note that  
$\wh{\alpha}=\alpha$ and hence $X^A$ is conformally invariant.

Given a choice of conformal scale,
the {\em tractor-$D$ operator}
$$
D_A\colon\ce_{B \cdots E}[w]\to\ce_{AB\cdots E}[w-1]
$$
is defined by 
\begin{equation}\label{Dform}
D_A V:=(n+2w-2)w Y_A V+ (n+2w-2)Z_{Aa}\nabla^a V -X_A\Box V, 
\end{equation} 
 where $\Box V :=\Delta V+w \J V$.  This also turns out to be
 conformally invariant as can be checked directly using the formulae
 above (or alternatively there are conformally invariant constructions
 of $D$, see e.g.\ \cite{Gosrni}).

The curvature $ \Omega$ of the tractor connection 
is defined by 
$$
[\nd_a,\nd_b] V^C= \Omega_{ab}{}^C{}_EV^E 
$$
for $ V^C\in\ce^C$.  Using
\eqref{basictrf} and the formulae for the Riemannian curvature yields
\begin{equation}\label{tractcurv}
\Omega_{abCE}= Z_C{}^cZ_E{}^e C_{abce}-2X_{[C}Z_{E]}{}^e A_{eab}
\end{equation}

\subsection{Forms and tensors}

The basic tractor tools for dealing with weighted differential forms
are developed in \cite{Branson+Gover} and following that source we
write $\ce^k[w]$ for the space of sections of $(\Lambda^k T^*M
)\otimes E[w]$ (and $\ce^k= \ce^k[0]$).  However in order to be
explicit and efficient in calculations involving bundles of possibly
high rank it is necessary to introduce some further abstract index
notation.  In the usual abstract index conventions one would write
$\ce_{[ab\cdots c]}$ (where there are implicitly $k$-indices skewed
over) for the space $\ce^k$. To simplify subsequent expressions we
 use the following conventions. Firstly 
indices labelled with sequential superscripts which are 
at the same level (i.e. all contravariant or all
covariant) will indicate a completely skew set of indices.  
Formally we set $a^1 \cdots a^k = [a^1 \cdots a^k]$ and so, for example,
$\ce_{a^1 \cdots a^k}$ is an alternative notation for $\ce^k$
while $\ce_{a^1 \cdots a^{k-1}}$ and $\ce_{a^2 \cdots a^k}$ both denote $\ce^{k-1}$. Next we
abbreviate this notation via multi-indices: We will use the forms
indices 
$$
\begin{aligned}
 \vec{a}^k &:=a^1 \cdots a^k =[a^1 \cdots a^k], \quad k \geq 0,\\
\dot{\vec{a}}^k &:= a^2 \cdots a^k=[a^2 \cdots a^k], \quad k \geq 1,\\
\ddot{\vec{a}}^k &:= a^3 \cdots a^k =[a^3 \cdots a^k], \quad k \geq 2,  \\
\dddot{\vec{a}}^k &:= a^4 \cdots a^k =[a^4 \cdots a^k], \quad k \geq 3.  \\
\end{aligned}
$$ If, for example, $k=1$ then $ \dot{\vec{a}}^k$ simply means the
index is absent, whereas if $k=1$ then $\ddform{a}$ means the term
containing the index $\ddform{a}$ is absent. For
example, a 3--form $\ph$ can have the following possible equivalent
structures of indices:
$$ \ph_{a^1a^2a^3} = \ph_{[a^1a^2a^3]} = \ph_{\vec{a}^3} 
= \ph_{a^1\dot{\vec{a}}^2} 
= \ph_{[a^1\dot{\vec{a}}^3]} = 
   \ph_{a^1a^2\ddot{\vec{a}}^3} \in \mathcal{E}_{\vec{a}^3} = \mathcal{E}^3. 
$$
We will also use $\bg_{\vec{a}^k\vec{b}^k}$
(and similarly $\bg_{\dot{\vec{a}}^k\dot{\vec{b}}^k}$) for $\bg_{a^1b^1}
\cdots \bg_{a^kb^k}$ (where all $a$--indices and all $b$--indices are
skewed over) and $\bg$ denotes the conformal metric.

The corresponding notations will be used for tractor indices so
e.g. the bundle of tractor $k$--forms $\ce_{[A^1\cdots A^k]}$ will be
denoted by $\ce_{A^1\cdots A^k}$ or $\mathcal{E}_{\vec{A}^k}$.

\vspace{1ex}

We shall demonstrate the notation by giving the conformal 
transformation formulae
of the Levi--Civita connection acting on conformally weighted forms. 
Under a 
rescaling $g \mapsto \wh{g} = e^{2\Up}g$ of the metric, and writing 
$\Up_a : = \na_a \Up$,  from (\ref{grad_trans_gen})
we have 
\begin{equation} \label{grad_trans_form}
\begin{split}
  \hat{\na}_{a^0} f_{\form{a}^k} &= \na_{a^0} f_{\form{a}^k}
    + w \Up_{a^0} f_{\form{a}^k} \\
  \hat{\na}^{a^1} f_{\form{a}^k} &= \na^{a^1} f_{\form{a}^k}
    + (n+w-2k) \Up^{a^1} f_{\form{a}^k},
\end{split}
\end{equation}
for $f_{\form{a}^k} \in \cE_{\form{a}^k}[w]$.

We will need similar results for spaces with more complicated 
symmetries. We shall define $\cE(1,k)$ for $k \geq 1$ and
$\cE(2,k)$ for $k \geq 2$ as follows:
\begin{align*}
  & \cE(1,k) := \{ f_{c\form{a}^k} \in  \cE_{c\form{a}^k} \mid 
      f_{[c\form{a}^k]}=0 \} \subseteq \cE_{c\form{a}^k} \\
  & \cE(2,k) := \{ \tilde{f}_{\form{c}^2\form{a}^k}  
      \in \cE_{\form{c}^2\form{a}^k} \mid 
      \tilde{f}_{[\form{c}^2\form{a}^k]} = 
      \tilde{f}_{c^1[c^2\form{a}^k]} =
      \tilde{f}_{[\form{c}^2\form{a}^{k-1}]a^k} =0\}
      \subseteq \cE_{\form{c}^2\form{a}^k}.
\end{align*}
In  other words, the subspaces $\cE(1,k)$ and $\cE(2,k)$ 
are defined by the condition that any skew symmetrisation of more than 
$k$ indices vanishes. The subspaces of completely trace-free tensors in 
$\cE(1,k)$ and $\cE(2,k)$ 
 will be denoted respectively by 
$\cE(1,k)_0$ and $\cE(2,k)_0$.  Tensor products with density bundles
will be denoted in an obvious way. For example $\cE(1,k)_0[w]$ is a shorthand for $\cE(1,k)_0 \otimes \ce[w]$. 

 We will later need the following identities
\begin{equation} \label{skew}
  f_{a^1p\dform{a}^k} = \frac{1}{k} f_{p\form{a}^k}
  \quad \mbox{and} \quad
  \tilde{f}_{a^1qp\dform{a}^k} = \frac{1}{k} \tilde{f}_{pq\form{a}^k}
\end{equation}
for $f_{c\form{a}^k} \in \cE(1,k)[w]$ and $\tilde{f}_{\form{c}^2
\form{a}^k} \in \cE(2,k)[w]$.  This follows from the skewing
$[p\form{a}^k]$ which vanishes in both cases. Using the second of
these we recover, for example, the well known identities
$$
R_{[a\ c]}^{\;\ b\;\ d} = 
  \frac{1}{2} R_{ac}^{\ \ bd} \quad \mbox{and} \quad C_{[a\ c]}^{\;\ b\;\ d} = 
  \frac{1}{2} C_{ac}^{\ \ bd} ~.
$$
Via \nn{skew}, (\ref{grad_trans_gen}) and  a short
computation we obtain the transformations 
\begin{equation} \label{grad_trans_(j,k)}
\begin{split}
  \hat{\na}_{a^0} f_{c\form{a}^k} 
  & = \na_{a^0} f_{c\form{a}^k}
    + (w-1) \Up_{a^0} f_{c\form{a}^k} 
    + \bg_{ca^0} \Up^p  f_{p\form{a}^k} \\
  \hat{\na}^{c} f_{c\form{a}^k}
  & =  \na^c f_{c\form{a}^k}
    + (n+w-k-1) \Up^c f_{c\form{a}^k} \\
  \hat{\na}^{c^1} \tilde{f}_{\form{c}^2 \form{a}^k} 
  & =  \na^{c^1} \tilde{f}_{\form{c}^2\form{a}^k}
   + (n+w-k-3) \Up^{c^1} \tilde{f}_{\form{c}^2\form{a}^k}
\end{split}
\end{equation}
for $f_{c\form{a}^k} \in \cE(1,k)_0[w]$ and 
$\tilde{f}_{\form{c}^2 \form{a}^k} \in \cE(2,k)_0[w]$.

\subsection{Tractor forms.} \label{s.form_tractor}
It follows 
from the semidirect composition series  of $\mathcal{E}_{A}$ that  the
corresponding decomposition of $\mathcal{E}_{\vec{A}^k}$ is
\begin{equation} \label{comp_series_form}
  \mathcal{E}_{[A^1 \cdots A^k]} = \mathcal{E}_{\vec{A}^k} \simeq
  \mathcal{E}^{k-1}[k] \lpl \left( \mathcal{E}^k[k] \oplus
  \mathcal{E}^{k-2}[k-2] \right) \lpl \mathcal{E}^{k-1}[k-2].
\end{equation}
 Given a choice of metric $g$ from the conformal class this
determines a splitting of this space into four components (a
replacement of the $\lpl$s with $\oplus$s is effected) and the
projectors (or splitting operators) $X,Y,Z$ for $\mathcal{E}_A$
determine corresponding projectors $\X,\Y,\Z,\W$ for
$\mathcal{E}_{\vec{A}^{k+1}}$, $k \geq 1$ as follows.
\begin{center}
\renewcommand{\arraystretch}{1.3}
\begin{tabular}{c@{\;=\;}l@{\;=\;}l@{\;=\;}l@{\ $\in$\ }l}
$\Y^k$ & $\Y_{A^0A^1 \cdots A^k}^{\quad a^1 \cdots\, a^k}$ &
  $\Y_{A^0\vec{A}^k}^{\quad \vec{a}^k}$ & $Y_{A^0}^{}Z_{A^1}^{a^1} \cdots Z_{A^k}^{a^k}$ &
  $\mathcal{E}_{\vec{A}^{k+1}}^{\vec{a}^k}[-k-1]$ \\
$\Z^k$ & $\Z_{A^1 \cdots A^k}^{\, a^1 \cdots\, a^k}$ &
  $\Z_{\vec{A}^k}^{\,\vec{a}^k}$ & $Z_{A^1}^{\,a^1} \cdots Z_{A^k}^{\,a^k}$ &
  $\mathcal{E}_{\vec{A}^k}^{\vec{a}^k}[-k]$ \\
$\W^k$ & $\W_{A'A^0A^1 \cdots A^k}^{\quad\,\ \ a^1 \cdots\, a^k}$ &
  $\W_{A'A^0\vec{A}^k}^{\quad\,\ \ \vec{a}^k}$ &
  $X_{[A'}^{}Y_{A^0}^{}Z_{A^1}^{\,a^1} \cdots Z_{A^k]}^{\,a^k}$ &
  $\mathcal{E}_{\vec{A}^{k+2}}^{\vec{a}^k}[-k]$ \\
$\X^k$ & $\X_{A^0A^1 \cdots A^k}^{\quad a^1 \cdots\, a^k}$ &
  $\X_{A^0\vec{A}^k}^{\quad \vec{a}^k}$ & $X_{A^0}^{}Z_{A^1}^{\,a^1} \cdots Z_{\,A^k}^{a^k}$ &  
  $\mathcal{E}_{\vec{A}^{k+1}}^{\vec{a}^k}[-k+1]$
\end{tabular}
\end{center}
where $k \geq 0$. The superscript $k$ in $\Y^k$, $\Z^k$, $\W^k$
and $\X^k$ shows always the corresponding tensor valence. (This is
slightly different than in \cite{Branson+Gover}, where $k$ concerns
the tractor valence.) Note that $Y=\Y^0$, $Z=\Z^1$ and $X=\X^0$ and
$\W^0 = X_{[A'}Y_{A^0]}$.
Using these projectors a section $f_{\vec{A}^{k+1}} \in
\mathcal{E}_{\vec{A}^{k+1}}$ can be written as a 4-tuple
$$f_{\vec{A}^{k+1}} =\begin{pmatrix}
    \si_{\form{a}^k} \\ \mu_{a^0 \form{a}^k}\ \ph_{\dform{a}^k} \\ \rh_{\form{a}^k}
  \end{pmatrix} =
  \Y_{A^0\vec{A}^k}^{\quad \vec{a}^k} \si_{\vec{a}^k} +
  \Z_{A^0\vec{A}^k}^{\,a_0\, \vec{a}^k} \mu_{a_0\vec{a}^k} +
  \W_{A^0\vec{A}^k}^{\quad \dot{\vec{a}}^k} \ph_{\dot{\vec{a}}^k} +
  \X_{A^0\vec{A}^k}^{\quad \vec{a}^k} \rh_{\vec{a}^k} $$
for forms $\si,\mu,\ph,\rh$ of weight and valence according to the relationship 
given in (\ref{comp_series_form}). 

The conformal transformation \nn{transf-tractor} yields the
transformation formulae for the projectors:
\begin{align} \label{transformation}
\begin{split}
  \widehat{\Y_{A^0\vec{A}^k}^{\quad \vec{a}^k}} 
    =&   \Y_{A^0\vec{A}^k}^{\quad \vec{a}^k}
       - \Up_{a^0}\Z_{A^0\vec{A}^k}^{\,a^0\,\vec{a}^k}
       - k \Up^{a^1}\W_{A^0\vec{A}^k}^{\quad \dot{\vec{a}}^k} \\
     & - \frac{1}{2} \Up^k\Up_k^{}\X_{A^0\vec{A}^k}^{\quad \vec{a}^k}
       + k \Up_p^{}\Up^{a^1} \X_{A^0\vec{A}^k}^{\quad \!p\dot{\vec{a}}^k} \\
  \widehat{\Z_{A^0\vec{A}^k}^{\,a^0\,\vec{a}^k}} 
    =&   \Z_{A^0\vec{A}^k}^{\,a^0\,\vec{a}^k}
       + (k+1) \Up^{a^0} \X_{A^0\vec{A}^k}^{\quad\, \vec{a}^k} \\
  \widehat{\W_{A^0\vec{A}^k}^{\quad \dot{\vec{a}}^k}} 
    =&   \W_{A^0\vec{A}^k}^{\quad \dot{\vec{a}}^k}
       - \Up_{a^1}^{} \X_{A^0\vec{A}^k}^{\quad \vec{a}^k} \\
  \widehat{\X_{A^0\vec{A}^k}^{\quad \vec{a}^k}} 
    =&   \X_{A^0\vec{A}^k}^{\quad \vec{a}^k}
\end{split}
\end{align}
for metrics $\hat{g}$ and $g$ from the conformal class. The normal tractor 
connection on $(k+1)$-form-tractors is
\begin{equation}
\label{normtrconn}
 \na_p \begin{pmatrix}
    \si_{\form{a}^k} \\ \mu_{a^0\form{a}^k}\ \ph_{\dform{a}^k} 
      \\ \rh_{\form{a}^k}
  \end{pmatrix} =
  \begin{pmatrix}
    \na_p \si_{\vec{a}^k} - (k+1) \mu_{p\vec{a}^k} - \bg_{pa^1} 
      \ph_{\dot{\vec{a}}^k} \\[2mm]
    \Bigl\{ \begin{smallmatrix}
               \na_p \mu_{a^0\vec{a}^k} \\ 
               + P_{pa^0} \si_{\vec{a}^k} + g_{pa^0} \rh_{\vec{a}^k}
            \end{smallmatrix} \Bigr\} \quad
    \Bigl\{ \begin{smallmatrix}
               \na_p \ph_{\dot{\vec{a}}^k} \\ 
               + k P_p^{a^1} \si_{\vec{a}^k} - k \de_p^{a^1} \rh_{\vec{a}^k}    
            \end{smallmatrix} \Bigr\} \\[2mm]
    \na_p \rh_{\vec{a}^k} - (k+1) P_p^{\ a^0} \mu_{a^0\vec{a}^k}^{} 
      + P_{pa^1} \ph_{\dot{\vec{a}}^k}  
  \end{pmatrix} 
\end{equation}
or equivalently 
\begin{eqnarray*}
  \na_p \Y_{A^0\vec{A}^k}^{\quad \vec{a}^k} &=&
      P_{pa_0} \Z_{A^0\vec{A}^k}^{\,a^0\,\vec{a}^k}
    + k P_p^{\ a^1} \W_{A^0\vec{A}^k}^{\quad \dot{\vec{a}}^k} \\
  \na_p \Z_{A^0\vec{A}^k}^{\,a^0\,\vec{a}^k} &=&
    - (k+1) \de_p^{a^0} \Y_{A^0\vec{A}^k}^{\quad \vec{a}^k}
    - (k+1) P_p^{\ a^0}  \X_{A^0\vec{A}^k}^{\quad \vec{a}^k} \\
  \na_p \W_{A^0\vec{A}^k}^{\quad\, \dot{\vec{a}}^k} &=&
    - \bg_{pa^1} \Y_{A^0\vec{A}^k}^{\quad \vec{a}^k}
    + P_{pa^1} \X_{A^0\vec{A}^k}^{\,\ \ a^1\!\dot{\vec{a}}^k} \\
  \na_p \X_{A^0\vec{A}^k}^{\quad \vec{a}^k} &=&
      \bg_{pa^0} \Z_{A^0\vec{A}^k}^{\,a^0\,\vec{a}^k}
    - k \de_p^{a^1} \W_{A^0\vec{A}^k}^{\quad \dot{\vec{a}}^k} .
\end{eqnarray*}

\subsection{The conformal Killing equation on forms.} 
The space $\cE_{c\vec{a}^k} = \ce_c  \otimes \ce_{a^1\cdots a^k}$ 
is completely reducible for $1 \leq k \leq n$ and we have the
$O(g)$-decomposition $\mathcal{E}_{c\vec{a}^k}[w] \cong
\mathcal{E}_{[c\vec{a}^k]}[w] \oplus \mathcal{E}_{\{c\vec{a}^k\}_0}[w]
\oplus \mathcal{E}_{\vec{a}^{k-1}}[w-2]$ 
where the bundle $ \mathcal{E}_{\{c\vec{a}^k\}_0}[w]$ consists of rank
$k+1$ trace-free tensors $T_{c\vec{a}^k}$ (of conformal weight $w$)
that are skew on the indices $a^1\cdots a^k$ and have the property
that $T_{[ca^1\cdots a^k]}=0$.  (Note that the three spaces on the
right-hand side are $SO(g)$-irreducible if $k \not \in \{ n/2, n/2 \pm
1 \}$). 
On the space $\mathcal{E}_{c\vec{a}^k}[w] $ there is a projection 
$\cP_{\{c\vec{a}^k\}_0}$ to the 
component  $\mathcal{E}_{\{c\vec{a}^k\}_0}[w]$ and we will use the notation 
$$
T_{c\vec{a}^k} \stackrel{\{c\vec{a}^k\}_0 }{=} S_{c\vec{a}^k} \quad \mbox{ or }
\quad T_{c\vec{a}^k} {=}_{\{c\vec{a}^k\}_0 } S_{c\vec{a}^k} 
$$
to mean that $\cP_{\{c\vec{a}^k\}_0}(T)=\cP_{\{c\vec{a}^k\}_0}(S)$. 
We will also use the projection $\cP_{\{c\vec{a}^k\}}$ to 
$\ce(1,k)[w]=:\cE_{\{c\vec{a}^k\}}[w]$.

Each metric from the conformal class determines a corresponding
Levi-Civita connection $\nd$ and for $1 \leq k \leq n-1$ and
$\si_{\vec{a}^k} \in \mathcal{E}^k[k+1]$, we may form $\nd_c
\si_{\vec{a}^k}$. This is not conformally invariant. However it is
straightforward to verify that its projection $\cP_{\{c\vec{a}^k
\}_0}(\na \si)$ is conformally invariant. That is, this  is
independent of the choice of metric (and corresponding Levi-Civita
connection) from the conformal class.  Thus the equation
\begin{equation}\label{ke}
  \na_{\{c} \si_{\vec{a}^k\}_0} =0, \quad \quad 1\leq k\leq n-1
  \tag{CKE}
\end{equation}
called the (form) \idx{conformal Killing equation}, is conformally
invariant. This is exactly the equation \nn{ifv} from the introduction. 

Suppose $\tilde{\nd}$ is a connection on another vector bundle (or
space of sections thereof) $\ce_\bullet$. For this connection coupled
with the Levi-Civita connection let us also write $\tilde{\nd}$. Since
it is a first order equation \nn{ke} is {\em strongly invariant} (cf.\
\cite{GoSrni04,MikESrni}) in the sense that if now $\si_{\vec{a}^k}
\in \mathcal{E}_{\vec{a}^k\bullet}[k+1]=
\mathcal{E}_{\vec{a}^k}[k+1]\otimes \ce_\bullet$ then
$\tilde{\na}_{\{c} \si_{\vec{a}^k\}_0} =0$ is also conformally
invariant. We will also call any such equation a conformal Killing
equation (or sometimes for emphasis a {\em coupled conformal Killing
equation}).

On oriented conformal manifolds the conformal Hodge-$\star$ operator
(see e.g. \cite{Branson+Gover}) gives a mapping
$$
\star: \ce^k[w] \to \ce^{n-k}[n-2k+w] \quad k\in \{0,1,\cdots ,n\} ~.
$$
In particular we have 
$$
\star: \ce^k[k+1] \to \ce^{n-k}[n-k+1] ,
$$
and from elementary classical SO$(n)$-representation theory it
follows easily that $ \si\in \ce^k[k+1] $ solves \nn{ke} for $k$-forms
if and only if $\star \si$ solves the version of \nn{ke} for
$(n-k)$-forms. Thus on oriented manifolds it is only strictly
necessary to study this equation for (weighted) $k$-forms with $k\leq
n/2$. Also it follows that on even dimensional oriented manifolds a
form in $\ce^{n/2}[n/2+1]$ is a solution of \nn{ke} if and only if its
self-dual and anti-self dual parts are separately solutions.
Nevertheless, since the redundancy does us no harm, we shall ignore
these observations and in the following simply treat the equation on
$k$-forms for $1 \leq k \leq n-1$.

\section{Invariant prolongation for conformal Killing forms}
\label{fullpro}

Throughout this section, and in much of the subsequent work, we will
write $f_{\form{a}}$ (rather than $f_{\form{a}^k}$) to denote a
section in $\cE_{\form{a}^k}[k+1]$.  That is, the superscript of the
form index $\form{a}$ will be omitted but can be taken to be $k$ (or 
otherwise if clear from the context). 

Before we start with the construction of the prolongation, we will
introduce some notation for certain algebraic actions of the curvature on 
tensors. Let us write $\sharp $ (which we will term {\em hash}) for the natural 
action of sections $A$ of ${\rm End}(TM)$ on tensors. For example,
on a covariant 2-tensor $T_{ab}$, we have
$$
A\sharp T_{ab}=-A^{c}{}_aT_{cb}-A^{c}{}_bT_{ac}.
$$
If $A$ is skew for a metric $g$, then at each point, $A$ is
$\frak{so}(g)$-valued.  The hash action thus commutes with the raising
and lowering of indices and preserves the ${\rm SO}(g)$-decomposition
of tensors. For example the Riemann tensor may be viewed as an ${\rm
End}(TM)$-valued 2-form $R_{ab}$ and in this notation we have
$$
[\nd_a, \nd_b]T=R_{ab}\sharp T ~,
$$ for an arbitrary tensor $T$. Similarly we have $C_{ab}\sharp T$ for
the Weyl curvature.  As a section of the tensor square of the $g$-skew
bundle endomorphisms of $TM$, the Weyl curvature also has a double
hash action that we denote $C\sharp\sharp T$.

We need some more involved actions of the Weyl tensor on
$\cE_{\form{a}^k}[w]$ for $k \geq 2$. These are given by
\begin{equation} \label{lozenge}
\begin{split}
  (C \blacklozenge f)_{c\dform{a}} 
  :=& \frac{k-2}{k} \left(
      C_{ca^2}^{\quad\  pq} f_{pq\ddform{a}}^{} 
      + C_{a^3a^2}^{\quad\  pq} f_{pqc\dddform{a}}^{} \right)
      \in \cE_{c\dform{a}^k}[w-2] \\
  (C \lozenge f)_{\form{ca}}
  :=& C_{c^1c^2a^1}^{\qquad p} f_{p\dform{a}}^{} 
      + C_{a^1a^2c^1}^{\qquad p} f_{pc^2\ddform{a}}^{} 
      + \frac{k}{n-k} \bg_{c^1a^1} (C \blacklozenge f)_{c^2\dform{a}}^{} \\
    & \in \cE_{\form{c}^2\form{a}^k}[w]
\end{split}
\end{equation}
where $\form{c} = \form{c}^2$ and $f_{\form{a}} \in
\cE_{\form{a}^k}[w]$. Note that $C \lozenge f$ vanishes for $k=n-1$
since $\cE(2,n-1)_0$ is trivial. For the sake of complete clarity we have
given these explicit formulae but note that, up to a multiple, the
first of these is simply $C\sharp f \in \cE_{\form{c}^2\form{a}^k}$
followed by projection to $\ce(1,k-1)[w-2]$ (the projection involves a
trace),
while the second is
$C\sharp f$ followed by projection to $\ce(2,k)_0[w]$.  This is
clear except for the final projection in each case which we now
verify.

\begin{lem} \label{l.lozenge}
Let us suppose $k \geq 2$. Then

\mbox{\rm (i)} $(C \blacklozenge f)_{c\dform{a}} =  
C_{\{ca^2}^{\quad\  pq} f_{|pq|\ddform{a}\}}^{} \in \cE(1,k-1)_0[w-2]$

\mbox{\rm (ii)} $(C \lozenge f)_{\form{ca}} \in \cE(2,k)_0[w]$
\end{lem}
\begin{proof}
(i) It follows from (\ref{lozenge}) and the Bianchi identity that
$(C \blacklozenge f)_{c\dform{a}}$ is trace-free. Moreover 
\begin{equation} \label{black}
  C_{\{ca^2}^{\quad\  pq} f_{|pq|\ddform{a}\}}^{}
  = C_{ca^2}^{\quad pq} f_{pq\ddform{a}}^{}
    - C_{[ca^2}^{\quad\ pq} f_{|pq|\ddform{a}]}^{}
  = (C \blacklozenge f)_{c\dform{a}}.
\end{equation}
where the first equality is just the definition of the projection $\{..\}$
and the second follows from re-expressing of the skew symmetrisation
$[c\dform{a}]$ in the last display.

(ii) According to the definition of $\cE(2,k)_0$, we are required to show
that $(C \lozenge f)_{c^1[c^2\form{a}]} = 
(C \lozenge f)_{[\form{c}\dform{a}]a^{k+1}} =0$
(note  $(C \lozenge f)_{[\form{ca}]}=0$ is obvious from
(\ref{lozenge})) and also that $C \lozenge f $ is trace-free.
Both skew symmetrisation's $[c^2\form{a}]$ and $[\form{c}\dform{a}]$
kill the last term of $C \lozenge f$ in (\ref{lozenge}),
because $(C \blacklozenge f)_{[c\dform{a}]}=0$ according
to the Lemma (i). Applying the symmetrisation $[c^2\form{a}]$ to the first two terms
in (\ref{lozenge}) and using the Bianchi identity yields
$$ C_{c^1[c^2a^1}^{\qquad p} f_{|p|\dform{a}]}^{}
   + \frac{1}{2} C_{[a^1a^2|c^1}^{\qquad\, p} f_{p|c^2\ddform{a}]}^{}, $$
where the indices $c^1c^2$ are {\em not} skewed over. This is zero because
$ C_{c^1[c^2a^1]}^{\qquad\, p} = -\frac{1}{2} C_{c^2a^1c^1}^{\qquad\, p}$.
The second skew symmetrisation $[\form{c}\dform{a}]$ is similar.

It remains to prove $\bg^{c^1a^1} (C \lozenge f)_{\form{ca}} =0$.
Tracing the last term in (\ref{lozenge}) yields
$$ \frac{k}{n-k} \bg^{c^1a^1} \bg_{c^1a^1} 
     (C \blacklozenge f)_{c^2\dform{a}}^{}
   = \frac{1}{2} (C \blacklozenge f)_{c^2\dform{a}}^{} $$
after a short computation. Further computations reveal
$$ \bg^{c^1a^1} C_{c^1c^2a^1}^{\qquad p} f_{p\dform{a}}^{} 
   =  -\frac{k-1}{2k} C_{c^2a^2}^{\quad\  pq} f_{pq\ddform{a}}^{} $$
and 
$$ \bg^{c^1a^1} C_{a^1a^2c^1}^{\qquad p} f_{pc^2\ddform{a}}^{}
   =  -\frac{k-2}{2k} C_{a^3a^2}^{\quad\  pq} f_{pqc^2\dddform{a}}^{} 
      +\frac{1}{2k} C_{c^2a^2}^{\quad\  pq} f_{pq\ddform{a}}^{}. $$
Summing the last three displays, the Lemma part (ii) follows from
(\ref{lozenge}) for $C \blacklozenge f$.
\end{proof}

\vspace{1ex}

Introducing new variables, the equation \nn{ke} may be
re-expressed in the form
$$
\na_c \si_{\form{a}}
  = \mu_{c\form{a}} + \bg_{ca^1} \nu_{\dform{a}} ~,
$$ where $\mu_{a^0\form{a}} \in \ce_{a^0\form{a}^k}[k+1]$ and
$\nu_{\dform{a}} \in \cE_{\dform{a}^k}[k-1]$. These capture some of
the 1-jet information: we have $\mu_{a^0\form{a}} = \na_{a^0}
\si_{\form{a}}$, and $\nu_{\dform{a}} = \frac{k}{n-k+1} \na^p
\si_{p\dform{a}}$. We need a further set of variables to complete 
\nn{ke} to a
first order closed system. There is some choice here, but, for the
purposes of studying conformal invariance, it turns out that 
$\rh_{\form{a}} := -\frac{1}{k}\na_{a^1} \nu_{\dform{a}} +
\frac{1}{nk} \na^p \na_{\{p} \si_{\form{a}\}_0} - P_{a^1}^{\;\
p} \si_{p\dform{a}}^{}$ is a judicious choice. We then have the
following result.

\begin{prop} \label{p.prolong_start}
Solutions of the conformal Killing equation \nn{ke}, for $1\leq k \leq n-1$,
are in 1-1 correspondence with solutions of the following system
on $\si_{\form{a}} \in \mathcal{E}_{\form{a}^k}[k+1]$,
$\mu_{a^0\form{a}} \in \ce_{a^0\form{a}^k}[k+1]$,
$\nu_{\dform{a}}  \in \ce_{\dform{a}^k}[k-1]$ and 
$\rh_{\form{a}}\in \ce_{\form{a}^k}[k-1]${\rm :}
\begin{equation} \label{prolong_start}
\begin{split}
  & \na_c \si_{\form{a}} = 
    \mu_{c\form{a}} + \bg_{ca^1} \nu_{\dform{a}} ~; \\
  & \na_c \mu_{a^0\form{a}} =
    (k+1) \left[ \bg_{ca^0} \rh_{\form{a}}^{} 
                 - P_{ca^0} \si_{\form{a}} 
                 - \frac{1}{2}C_{a^0a^1c}^{\qquad p} \si_{p\dform{a}}^{} 
          \right] ~; \\
  & \na_c \nu_{\dform{a}} = 
    - k \left[ \rh_{c\dform{a}} + P_c^{\ p} \si_{p\dform{a}}^{} \right] 
    + \frac{k(k-1)}{2(n-k)} (C \blacklozenge \si)_{c\dform{a}} ~; \\
  & \na_c \rh_{\form{a}} =
    P_{ca^1} \nu_{\dform{a}} - P_c^{\ p} \mu_{p\form{a}}^{} + 
    \frac{1}{2} A^p_{\ a^1a^2} \si_{pc\ddform{a}}^{} 
    - A^p_{\ ca^1} \si_{p\dform{a}}^{} \\
  & \qquad\quad\ +\frac{1}{2} C_{a^1a^2c}^{\qquad p} \nu_{p\ddform{a}}^{}
    - \frac{k}{2(n-k)} \na_{a^1}^{} (C \blacklozenge \si)_{c\dform{a}} 
    \qquad \mbox{\rm for}\ k \geq 2; \\
  & \na_c \rh_{a^1} = P_{ca^1} \nu - P_c^{\ p} \mu_{pa^1}^{}
    + A_{a^1pc} \si^p
    \qquad\qquad\qquad\quad \mbox{\rm for}\ k=1.
\end{split}
\end{equation} 
The mapping from solutions $\si_{\form{a}}$ of \nn{ke} to
solutions $(\si_{\form{a}} , \mu_{a^0\form{a}}, \nu_{\dform{a}},
\rh_{\form{a}})$ of the system above is
\begin{equation} \label{promap}
\begin{split}
  \si_{\form{a}} \mapsto \Bigl(
  & \si_{\form{a}} ,~ \na_{a^0}\si_{\form{a}},~ 
    \frac{k}{n-k+1} \na^p\si_{p\dform{a}}, \\
  & \frac{1}{n k} \na^p \na_{\{p} \si_{\form{a}\}_0} 
      - \frac{1}{n-k+1} \na_{a^1} \na^p \si_{p\dform{a}} 
      - P_{a^1}^{\ p} \si_{p\dform{a}}^{} \Bigr) 
\end{split}
\end{equation}
\end{prop}

\begin{proof}
As mentioned above the first equation $\na_c \si_{\form{a}} =
\mu_{c\form{a}} + \bg_{ca^1} \nu_{\dform{a}}$ is simply a restatement
of the conformal Killing equation (\ref{ke}) afforded by introducing
the new variables $\mu_{a^0\form{a}} \in \ce_{[a^0\form{a}]}[k+1]$ and
$\nu_{\dform{a}}\in \ce_{\dform{a}}[k-1]$. (At this point and until
further notice below we take the rank of $\si$ to be in the range
$1\leq k\leq n-1$.) 

This equation also gives $\mu_{a^0\form{a}}$ and $\nu_{\dform{a}}$ in
terms of derivatives of $\si_{\form{a}}$. Thus the Proposition is
clear except that we should verify that if $\si_{\form{a}}$ solves
\nn{ke} then we have the second, third and fourth equations of
(\ref{prolong_start}).

To establish the second equation, let us observe
$(k+2) \na_{[c} \na_{a^0} \si_{\form{a}]} = \na_c \na_{a^0} \si_{\form{a}}
- (k+1) \na_{a^1} \na_{[a^0} \si_{c\dform{a}]}$, and that the left-hand-side 
vanishes due to the Bianchi identity. The first term on the right hand 
side is just $\na_c \mu_{a^0\form{a}}$ thus
\begin{align*} 
  \na_c \mu_{a^0\form{a}}      
  &= (k+1) \na_{a^1} \mu_{a^0c\dform{a}} 
   = (k+1) \na_{a^1} \left( \na_{a^0} \si_{c\dform{a}} 
           - \bg_{a^0[c} \nu_{\dform{a}]} \right) \\
  &= (k+1) \left( \frac{1}{2} R_{a^1a^0c}^{\qquad p} \si_{p\dform{a}}^{}
           - \frac{1}{k} \bg_{ca^0} \na_{a^1} \nu_{\dform{a}} \right)
\end{align*} 
where the second equality follows from the first equation in 
(\ref{prolong_start}) and the third equality from the Bianchi identity.
Now the equation for $\na_c \mu_{a^0\form{a}}$ in (\ref{prolong_start})
follows from the last display using (\ref{csplit}) and the relation
$\rh_{\form{a}} = - \frac{1}{k} \na_{a^1} \nu_{\dform{a}}
- P_{a^1}{}^{p}\si_{p\dform{a}}^{}$, which we have for solutions.

The second equation in (\ref{prolong_start}) concerns $\na_c
\nu_{\dform{a}} = \frac{k}{n-k+1} \na_c \na^p \si_{p\dform{a}}$.
Commuting the covariant derivatives we get $\na_c \na^p = R_c{}^{p}
\sharp + \na^p \na_c$ where, recall, $\sharp$ captures the action of
the Riemann curvature tensor $R$. Therefore
\begin{align*}
  (n\!-\!k\!+\!1) & \na_c \nu_{\dform{a}}
   = k \left[ R_{c\ p}^{\ p\ q}\si_{q\dform{a}}^{}
       + (k-1) R_{c\ a^2}^{\ p\,\ q} \si_{pq\ddform{a}}^{}  
       + \na^p \left( \mu_{cp\dform{a}} + \bg_{c[p} \nu_{\dform{a}]}
       \right) \right] \\
  &= k \bigl[ -\Ric_c^{\ p} \si_{p\dform{a}}^{} 
       + \frac{1}{2} (k-1) R_{ca^2}^{\quad pq} \si_{pq\ddform{a}}^{}
       - \na^p \mu_{pc\dform{a}}
       + \frac{1}{k} \na_c \nu_{\dform{a}} \bigr]
\end{align*}
where 
we have used $\na^p \nu_{p\ddform{a}} = \frac{k}{n-k+1} \na^p \na^q
\si_{qp\ddform{a}}=0$. Note that the last term here is a multiple of the 
left-hand-side. We consider the other terms on the right-hand-side. 
Recall that  (\ref{csplit}) gives $\Ric_{ab}=(n-2)\Rho_{ab}+\J\bg_{ab}$. 
Using (\ref{csplit}) also for the second term
on the right-hand-side, and the equation for $\na_c\mu_{a^0\form{a}}$
in (\ref{prolong_start}) for the third, a computation yields
\begin{align*}
  & - \Ric_c^{\ p} \si_{p  \dform{a}}^{} = - (n-2) P_c^{\ p} \si_{p\dform{a}}^{}
    - J \si_{c\dform{a}} \\
  & \frac{1}{2}(k-1) R_{ca^2}^{\quad pq} \si_{pq\ddform{a}}^{} = 
    \frac{1}{2}(k-1) C_{ca^2}^{\quad pq} \si_{pq\ddform{a}}^{}
    + 2 (k-1) \de_{\ [c}^p P_{a^2]}^{\;\ q} \si_{pq\ddform{a}}^{} \\
  & \qquad\qquad\quad = 
    \frac{1}{2}(k-1) C_{ca^2}^{\quad pq} \si_{pq\ddform{a}}^{}
    - (k-1) \left( P_{a^2}^{\;\ p} \si_{pc\ddform{a}}^{} 
                   - P_c^{\ p} \si_{p\dform{a}} \right) \\
  & - \na^p \mu_{pc\dform{a}} = - (n-k) \rh_{c\dform{a}} 
    + J \si_{c\dform{a}} - k P_{[c}^{\;\ p} \si_{|p|\dform{a}]}^{}
    - \frac{1}{2} (k-1) C_{[a^2c}^{\quad\ qp} \si_{|pq|\ddform{a}]}^{}. 
\end{align*}
 Hence the last but one display says that 
$\frac{n-k}{k} \na_c \nu_{\dform{a}}$
is equal to the sum of the right hand sides of the last display.
Now using the relation $- k P_{[c}^{\;\ p} \si_{|p|\dform{a}]}^{}
= -P_c^{\ p} \si_{p\dform{a}}^{} + 
(k-1) P_{a^2}^{\;\ p} \si_{pc\ddform{a}}^{}$ and 
(\ref{black}) we obtain immediately the third equation in 
(\ref{prolong_start}).

The last equation requires more computation. Let us first make an
observation about its skew-symmetric part $\na_{[c}
\rh_{\form{a}]}$. Using the definition of $\rh$ and the Bianchi
identity, we have $\na_{[c} \rh_{\form{a}]} = -\na_{[c}^{}
P_{a^1}^{\;\ p} \si_{|p|\dform{a}]}^{}$.  Using the Leibniz rule and
the first equation in (\ref{prolong_start}) for the right hand side, we
obtain
\begin{equation} \label{grad_rh_skew}
  \na_{[c} \rh_{\form{a}]} = 
    - \frac{1}{2} A^p_{\ [ca^1} \si_{|p|\dform{a}]}^{}
    - P_{[c}^{\;\ p} \mu_{|p|\form{a}]}^{},
\end{equation}
since the term $P_{a^1}^{\;\ p} \bg_{c[p}^{} \nu_{\dform{a}]}^{}$
vanishes after the skew symmetrisation $[c\form{a}]$. Now to compute
the full section $\na_c \rh_{\form{a}}$, we shall start with the 
equation for $\na_c \nu_{\dform{a}}$ from (\ref{prolong_start}). 
We apply $\na_{a^1}$ to both sides of this equation and skew over all
$a$--indices. Commuting the covariant derivatives on the left-hand-side,
we obtain $\na_{a^1} \na_c = \na_c \na_{a^1} + R_{a^1c} \sharp$.
The first term on the right hand side is $-k \na_{a^1} \rh_{c\dform{a}}
= (k+1) \na_{[c} \rh_{\form{a}]} - \na_c \rh_{\form{a}}$.
Through these observations, and using (\ref{grad_rh_skew}), we obtain
\begin{gather*}
  \na_c \na_{a^1} \nu_{\dform{a}} 
    + (k-1) R_{a^1ca^2}^{\qquad p} \nu_{p\ddform{a}}^{} = 
    - (k+1) \bigl( \frac{1}{2} A^p_{\ [ca^1} \si_{|p|\dform{a}]}^{}
                    + P_{[c}^{\;\ p} \mu_{|p|\form{a}]}^{} \bigr) \\
  - \na_c \rh_{\form{a}} 
     - k \na_{a^1}^{} P_c^{\ p} \si_{p\dform{a}}^{}
     + \frac{k(k-1)}{2(n-k)} \na_{a^1} (C \blacklozenge \si)_{c\dform{a}}.
\end{gather*} 
Many terms can be simplified and we shall start with the the first term on 
the left-hand-side. We have 
$$ \na_c \na_{a^1} \nu_{\dform{a}} =
   - k \bigl( \na_c \rh_{\form{a}}
   + \na_c^{} P_{a^1}^{\ p} \si_{|p|\dform{a}}^{} \bigr)  $$
which follows from the equation for $\na_c \nu_{\dform{a}}$
in (\ref{prolong_start}). Combining the last two displays we obtain
\begin{gather*}
  - (k-1) \na_c \rh_{\form{a}} = 
      2k \na_{[c}^{} P_{a^1]}^{\ \ p} \si_{p\dform{a}}^{}
    - (k+1) \bigl( \frac{1}{2} A^p_{\ [ca^1} \si_{|p|\dform{a}]}^{}
                    + P_{[c}^{\;\ p} \mu_{|p|\form{a}]}^{} \bigr) \\
    - \frac{1}{2} (k-1) R_{a^1a^2c}^{\qquad p} \nu_{p\ddform{a}}^{}
    + \frac{k(k-1)}{2(n-k)} \na_{a^1} (C \blacklozenge \si)_{c\dform{a}}.
\end{gather*}
where we have also used $ R_{a^1ca^2}^{\qquad p} = \frac{1}{2}
R_{a^1a^2c}^{\qquad p}$. Note that for the case of (the rank of $\si$
being) $k=1$ both sides of the equality above vanish and we get no
information.
Now we  simplify terms on the right hand side:
the first term using the Leibniz rule and the equation for 
$\na_c \si_{\form{a}}$, the next two terms re-expressing the skew 
symmetrisation $[c\form{a}]$ and the first curvature term
using the decomposition (\ref{csplit}). This yields
\begin{align*}
  & 2k \na_{[c}^{} P_{a^1]}^{\ \ p} \si_{p\dform{a}}^{} =
      k A^p_{\ ca^1} \si_{p\dform{a}}^{}
      + 2k P_{[a^1}^{\ \ p} \mu_{c]p\dform{a}}
      + 2k P_{[a^1}^{\ \ p} \bg_{c][p}^{} \nu_{\dform{a}]}^{} \\
  & \qquad\qquad = k A^p_{\ ca^1} \si_{p\dform{a}}^{}
      + k P_{a^1}^{\;\ p} \mu_{cp\dform{a}}^{}
      - k P_{c}^{\ p} \mu_{a^1p\dform{a}}^{}
      + (k-1) \bg_{ca^1}^{} P_{a^2}^{\;\ p} \nu_{p\ddform{a}}^{} \\
  & - \frac{1}{2} (k\!+\!1) A^p_{\ [ca^1} \si_{|p|\dform{a}]}^{} =
      - A^p_{\ ca^1} \si_{p\dform{a}}^{} 
      + \frac{1}{2} (k\!-\!1) A^p_{\ a^2a^1} \si_{pc\ddform{a}}^{} \\
  & - (k\!+\!1) P_{[c}^{\;\ p} \mu_{|p|\form{a}]}^{} = 
      - P_{c}^{\;\ p} \mu_{p\form{a}}^{}
      + k P_{a^1}^{\;\ p} \mu_{pc\dform{a}}^{} \\
  & - \frac{1}{2} (k\!-\!1) R_{a^1a^2c}^{\qquad p} \nu_{p\ddform{a}}^{} = 
      - \frac{1}{2} (k\!-\!1) \bigl[
      C_{a^1a^2c}^{\qquad p} \nu_{p\ddform{a}}^{}
      + 2 \bg_{ca^1}^{} P_{a^2}^{\;\ p} \nu_{p\ddform{a}}^{}
      + 2 P_{ca^1} \nu_{\dform{a}} \bigr].  
\end{align*}
Substituting these in the previous display, the Proposition  for $k \geq 2$ 
follows. The case $k=1$ 
 can be checked directly by tracing 
$ \frac{1}{2} R_{c^0c^1}\sharp \mu_{a^0a^1} = 
\na_{c^0} \na_{c^1} \mu_{a^0a^1} =
\na_{c^0} \left[ 2\bg_{c^1a^0} \rh_{a^1} -2 P_{c^1a^0} \si_{a^1} 
- C_{a^0a^1c^1}^{\qquad p}\si_p^{} \right].$ 
\end{proof}

\noindent {\bf Remark:} There is a variant of the derivation for the
$k\geq 2$ cases, as in the proof above, which generalises the
treatment of $k=1$ that we give there. However this breaks down for
$k=n-1$. Dually the proof we give for $k \geq 2$ breaks down at 
$k=1$.  Our proof of
the $k=1$ agrees with a treatment of that case distributed privately
by Mike Eastwood during the preparation of \cite{BCEG} and his
notation and conventions influenced our treatment. Earlier alternative
treatments of that case have been known to the first author for some
time (see \cite{powerslap}). 

\vspace{1ex}

\begin{lem} \label{l.whitelozenge}
Let us fix $k \geq 2$. If $\si_{\form{a}} \in \cE_{\form{a}^k}[k+1]$
is a solution of (\ref{ke}) then $(C \lozenge \si)_{\form{ca}}=0$.
\end{lem}
\begin{proof}
We shall prove the lemma using the prolongation
(\ref{prolong_start}). Applying $\na_{c^1}$ to both sides of the equation 
for $\na_{c^2} \si_{\form{a}}$, we obtain
$$ \na_{c^1} \na_{c^2} \si_{\form{a}} 
   = \na_{c^1} \mu_{c^2\form{a}} 
   + \bg_{c^2a^1} \na_{c^1} \nu_{\dform{a}}. $$
The left-hand side is equal to
$$ \frac{k}{2} R_{c^1c^2a^1}^{\qquad p} \si_{p\dform{a}}^{}
   = \frac{k}{2} C_{c^1c^2a^1}^{\qquad p} \si_{p\dform{a}}^{} 
     + k \bg_{c^1a^1}^{} P_{c^2}^{\ p} \si_{p\dform{a}}^{}
     + k P_{c^1a^1} \si_{c^2\dform{a}} $$ 
according to (\ref{csplit}). On the other hand, from 
(\ref{prolong_start}) the right-hand side is equal to 
\begin{gather*} 
  \Bigl( -k \bg_{c^1a^1} \rh_{c^2\dform{a}} 
    + k P_{c^1a^1} \si_{c^1\dform{a}}
    - \frac{1}{2} \bigl( 2 C_{c^2a^1c^1}^{\qquad p} \si_{p\dform{a}}^{}
        - (k-1)  C_{a^2a^1c^1}^{\qquad p} \si_{pc^2\ddform{a}}^{} 
      \bigr)
    \Bigr) \\
  + \bg_{c^2a^1} \Bigl( -k \rh_{c^1\dform{a}}
    - k P_{c^1}^{\ p} \si_{p\dform{a}}^{}
    + \frac{k(k-1)}{2(n-k)} (C \blacklozenge \si)_{c^1\dform{a}}
    \Bigr). 
\end{gather*} 
Now equating these two displays and 
using $C_{c^2a^1c^1}^{\qquad p} = -\frac{1}{2} C_{c^1c^2a^1}^{\qquad p}$
we obtain an identity which holds for solutions. Comparing the expression 
with  the definition 
of $(C \lozenge \si)$ in (\ref{lozenge}), we see the identity is 
$$
(k-1)(C \lozenge \si)=0. \vspace{-5.4ex}
$$
\end{proof}

Note that a curvature condition, equivalent to that in Lemma
\ref{l.whitelozenge}, is in \cite{Kashiwada}. There the identity for
solutions is stated in terms of the Riemann tensor $R$, rather than in
terms of the Weyl tensor $C$. In this form it has also been derived in
\cite{Sem-hab} (although we could not find the necessary restriction
$k \geq 2$ in that source).  Expressing the
identity via the Weyl curvature, as we do, emphasises that this is a
conformally invariant condition.

Next we observe that \nn{promap} defines a conformally invariant
differential splitting operator. We define a differential operator 
$ {\Bbb D}$ on  $\cE_{\form{a}^k}[k+1]$ by 
\begin{equation} \label{si_A}
  \si_{\form{a}} \mapsto \si_{A^0\form{A}} := 
    \Y_{A^0\form{A}}^{\quad \form{a}} \si_{\form{a}}^{} 
  + \frac{1}{k+1} \Z_{A^0\form{A}}^{\,a^0\,\form{a}} \mu_{a^0\form{a}}^{} 
  + \W_{A^0A^1\dform{A}}^{\quad\;\ \  \dform{a}} \nu_{\dform{a}}^{}  
  - \X_{A^0\form{A}}^{\quad \form{a}} \rh_{\form{a}}^{},
\end{equation}
where $\si_{\form{a}}$, $\mu_{a^0\form{a}}$, $\nu_{\dform{a}}$ and 
$\rh_{\form{a}}$ are given by \nn{promap}. Then we have the following.

\begin{lem} \label{D_invariance}
  ${\Bbb D}$ is a conformally invariant operator 
$$ {\Bbb D}: \cE_{\form{a}^k}[k+1] \to \cE_{A^0\form{A}^k} ~ \quad \mbox{for}\quad 1\leq k\leq n-1.$$
\end{lem} 

\begin{proof}
Consider $\D$ for $\si \in \cE_{\form{a}^k}[k+1]$. Let $\mu$, $\nu$
and $\rh$ be given in terms of $\si$ as in \nn{promap}. 
In these
formulae $\na$ is the Levi-Civita connection for some choice of metric
$g$ from the conformal class.  So $\mu$, $\nu$ and $\rh$ depend on the
metric. If we conformally rescale the metric $g \mapsto
\wh{g}=e^{2\Up}g $ then it is easy to calculate (using e.g. the
transformation formulae given in \cite{GoPetLap}) that the sections
$\wh{\mu}$ and $\wh{\nu}$ for the metric $\wh{g}$ are given by
$\wh{\mu}_{a^0\form{a}} = \mu_{a^0\form{a}} + (k+1)
\Up_{a^0}\si_{\form{a}}$ and $\wh{\nu}_{\dform{a}} = \nu_{\dform{a}} +
k \Up^p\si_{p\dform{a}}$, where $\Up_a = \na_a \Up$. To compute
$\wh{\rh}_{\form{a}} = - \frac{1}{k} \wh{\na}_{a^1}
\wh{\nu}_{\dform{a}} - \wh{P}_{a^1}^{\;\ p} \si_{p\dform{a}}^{} +
\frac{1}{nk} \wh{\na}^p \wh{\na}_{\{p} \si_{\form{a}\}_0}$ we use the
transformations
\begin{eqnarray*}
  \wh{\na}_{a^1} \wh{\nu}_{\dform{a}} &=& \wh{\na}_{a^1}
    (\nu_{\dform{a}} + k \Up^p \si_{p\dform{a}}) =
    (\na_{a^1}+(k-1)\Up_{a^1})(\nu_{\dform{a}} + k \Up^p
    \si_{p\dform{a}}) \\
 &=& \na_{a^1} \nu_{\dform{a}}+
    (k-1)\Up_{a^1}\nu_{\dform{a}} + k(\na_{a^1}\Up^p) \si_{p\dform{a}}
    +k\Up^p\na_{a^1}\si_{p\dform{a}}\\
&&+ k(k-1)\Up_{a^1}\Up^p
    \si_{p\dform{a}}\\ 
\wh{P}_{a^1}^{\;\ p}
    \si_{p\dform{a}}^{} &=& P_{a^1}^{\;\ p} \si_{p\dform{a}}^{} -
    (\na_{a^1} \Up^p) \si_{p\dform{a}} + \Up_{a^1} \Up^p
    \si_{p\dform{a}} - \frac{1}{2} \Up^p \Up_p \si_{\form{a}} \\
\wh{\na}^p \wh{\na}_{\{p} \si_{\form{a}\}_0} &=& \wh{\na}^p
    \na_{\{p} \si_{\form{a}\}_0} = \na^p \na_{\{p} \si_{\form{a}\}_0}
    + n \Up^p \na_{\{p} \si_{\form{a}\}_0}.
\end{eqnarray*}
See (\ref{grad_trans_form}) for the first of these, \nn{Rhotrans} for the second
 and (\ref{grad_trans_(j,k)}) for the last.
Summing the right-hand sides with the required coefficients 
(from \nn{promap}) we get,
$$ \wh{\rh}_{\form{a}} = \rh_{\form{a}} - \frac{k-1}{k}\Up_{a^1}\nu_{\dform{a}}
     - \Up^p \na_{a^1} \si_{p\dform{a}}
     - k\Up_{a^1} \Up^p \si_{p\dform{a}} 
     + \frac{1}{2} \Up^p \Up_p \si_{\form{a}}
     + \frac{1}{k} \Up^p \na_{\{p} \si_{\form{a}\}_0}. $$
Recall $\frac{1}{k} \Up^p \na_{\{p} \si_{\form{a}\}_0} = 
\Up^p \na_{\{a^1} \si_{p\dform{a}\}_0}$ using (\ref{skew})
therefore
$-\Up^p \na_{a^1} \si_{p\dform{a}} + 
\frac{1}{k} \Up^p \na_{\{p} \si_{\form{a}\}_0} 
= -\Up^p (\mu_{a^1p\dform{a}} + \bg_{a^1[p} \nu_{\dform{a}]})$.
From this and the previous display we obtain
$$ \wh{\rh}_{\form{a}} = \rh_{\form{a}}
     + \Up^p \mu_{p\form{a}}
     - \Up_{a^1} \nu_{\dform{a}}
     + \frac{1}{2} \Up^p \Up_p \si_{\form{a}}
     - k \Up_{a^1} \Up^p \si_{p\dform{a}}.  $$
Using this and the transformation properties
from (\ref{transformation}), a short computation shows that
${\Bbb D} (\si)$ is a section of $\cE_{A^0\form{A}^k}$
that does
not depend on the choice of the metric from the conformal class.
\end{proof}

\noindent {\bf Remarks:} 1. For $k=1$, ${\Bbb D}$ is just the
$w=1$ and special case of the operator ${\Bbb D}^{\beta a}$ from
section 5.1 of \cite{BrGoOpava}.  

2. Note that the operator $\Bbb D$ is not unique as an invariant
differential operator ``putting'' $\si_{\form{a}} \in
\cE_{\form{a}^k}[k+1]$ into the top slot of $F_{A^0\form{A}} \in
\cE_{A^0\form{A}^k}$ (i.e.\ a differential splitting operator with
left inverse $F_{A^0\form{A}} \mapsto (k+1)\X_{\quad
  \form{a}}^{A^0\form{A}} F_{A^0\form{A}}^{}$).  $\Bbb D$ can be
obviously modified by any multiple of $\X_{A^0\form{A}}^{\quad
  \form{a}} C_{a^1a^2}^{\quad\ pq} \si_{pq\ddform{a}}^{}$.  

\vspace{1ex}

Assume $k \geq 2$. 
We define a 1st order differential operator
$$ \Ph_c: \cE_{A^0\form{A}^k} \longrightarrow \cE_{cA^0\form{A}^k}
$$ that will turn up in our later calculations. Given a section
$F_{A^0\form{A}} \in \cE_{[A^0\form{A}^k]}$ which, for a choice $g \in
[g]$ of the metric in the conformal class, is convenient to take to
be in the form
\begin{equation} \label{F_A^0A}
  F_{A^0\form{A}} =
    \Y_{A^0\form{A}}^{\quad \form{a}} \si_{\form{a}}^{}
  + \frac{1}{k+1} \Z_{A^0\form{A}}^{\,a^0\,\form{a}} \mu_{a^0\form{a}}^{}
  + \W_{A^0A^1\dform{A}}^{\quad\;\ \  \dform{a}} \nu_{\dform{a}}^{}
  - \X_{A^0\form{A}}^{\quad \form{a}} \rh_{\form{a}}^{}, 
\end{equation}
we set
\begin{equation} \label{Ph}
\begin{split}
  \Ph_c(F_{A^0\form{A}}) :=
  & - \frac{1}{2} \Z_{A^0\form{A}}^{\,a^0\,\form{a}}
      C_{a^0a^1c}^{\qquad p} \si_{p\dform{a}}^{}
    + \frac{k(k-1)}{2(n-k)} \W_{A^0A^1\dform{A}}^{\quad\;\ \  \dform{a}}
      (C \blacklozenge \si)_{c\dform{a}} \\
  & + \X_{A^0\form{A}}^{\quad \form{a}} \Bigl[
        A^p_{\;\ ca^1} \si_{p\dform{a}}^{}
      - \frac{1}{2} A^p_{\;\ a^1a^2} \si_{pc\ddform{a}}^{} 
      - \frac{1}{2} C_{a^1a^2c}^{\qquad p} \nu_{p\ddform{a}}^{} \\
  &   \qquad\qquad 
      + \frac{k}{2(n-k)} \na_{a^1}  (C \blacklozenge \si)_{c\dform{a}}
    \Bigr] ~.
\end{split}
\end{equation}

Our aim is to construct a connection ${}^k\na$ on
$\cE_{A^0\form{A}^k}$ such that solutions $\si_{\form{a}}$ of
(\ref{ke}) correspond to sections of $\cE_{A^0\form{A}^k}$ that are
parallel according to ${}^k\na$. Let us start with the normal
tractor connection $\na$. Using the previous proposition, it is a
short and straightforward calculation to show that if
$\si_{\form{a}}$ is a solution of (\ref{ke}), $k \geq 2$ then 
$\na_c \D(\si)_{A^0\form{A}} = \Ph_c(\D(\si)_{A^0\form{A}}).$ 
Also, it is easy to verify (or see \cite{powerslap}) that for $k=1$, 
if $\si_{a^1}$ is a solution 
of (\ref{ke}) then $\na_c \D(\si)_{A^0A^1} = \Om_{pcA^0A^1} \si^p.$
This leads us to the following.

\begin{lem} \label{l.nonclosed_noninv}
{\rm (i)} Given a metric $g$ from the conformal class, the mapping
$$ \si_{\form{a}} \mapsto {\D}(\si)_{A^0\form{A}}, 
   \quad \mbox{with inverse} \quad
   F_{A^0\form{A}} \mapsto (k+1)\X_{\quad \form{a}}^{A^0\form{A}} 
     F_{A^0\form{A}}^{}~, $$
gives a 
%conformally invariant 
bijective mapping
between sections of $\si_{\form{a}} \in \cE_{\form{A}^k}[k+1]$
satisfying (\ref{ke}) and sections $F_{A^0\form{A}} \in
\cE_{A^0\form{A}^k}$ satisfying, 
\begin{alignat*}{2}
  &\na_c F_{A^0\form{A}} = \Ph_c(F_{A^0\form{A}})    & \qquad k \geq 2, \\
  &\na_c F_{A^0A^1} = \Om_{pcA^0A^1} \si^p           & \qquad k=1.
\end{alignat*}

{\rm (ii)} 
Upon a change 
of the metric $g \mapsto \wh{g} = e^{2\Up}g$, $\Phi_c$ transforms according to
$$ \wh{\Ph}_c(F_{A^0\form{A}}) = \Ph_c(F_{A^0\form{A}})  
   - \X_{A^0\form{A}}^{\quad \form{a}} 
   \Up^p (C \lozenge \si)_{pc\form{a}} $$
where $\Up_a = \na_a \Up$ and 
$\si_{\form{a}} = (k+1)\X_{\quad \form{a}}^{A^0\form{A}} 
F_{A^0\form{A}}^{}$.
\end{lem}

\begin{proof}
%The mapping $\si_{\form{a}} \mapsto {\D}(\si)_{A^0\form{A}}$ on
%$\cE_{\form{a}}[k+1]$ is well-defined (on conformal manifolds)
%as shown in Lemma \ref{D_invariance}. 
We have already observed that $\na_c \D(\si)_{A^0\form{A}} =
\Ph_c(\D(\si)_{A^0\form{A}})$ for solutions $\si$ of (\ref{ke}) for
$k \geq 2$, and the also the corresponding statement for $k=1$. 
On the other hand, looking at the coefficients of $\Y$ on both sides of
$\na_c F_{A^0\form{A}} = \Ph_c(F_{A^0\form{A}})$ we see this relation
implies that the ``top slot'' $\si_{\form{a}}:= (k+1)\X_{\quad
\form{a}}^{A^0\form{A}} F_{A^0\form{A}}^{}$ of $F$ is a solution of
(\ref{ke}). Thus the claimed bijective correspondence follows.

It remains to prove (ii). Let us consider a section $F_{A^0\form{A}}$
of the form (\ref{F_A^0A}) and a conformal rescaling
$g \mapsto \wh{g}$ as above. 
Collecting together the conformal transformation formulae for all the
relevant objects we have:
\begin{align} \label{trans_proof}
\begin{split}
  \wh{\mu}_{a^\form{a}} 
  =& \mu_{a^\form{a}} + (k+1) \Up_{a^0} \si_{\form{a}} \\
  \wh{\nu}_{\dform{a}}  
  =& \nu_{\dform{a}} + k \Up^p \si_{p\dform{a}} \\
  \wh{\Z}_{A^0\vec{A}}^{\,a^0\,\vec{a}} 
  =& \Z_{A^0\vec{A}}^{\,a^0\,\vec{a}}
     + (k+1) \Up^{a^0} \X_{A^0\vec{A}}^{\quad\, \vec{a}} \\
  \wh{\W}_{A^0A^1\dform{A}}^{\qquad\! \dform{a}} 
  =& \W_{A^0A^1\dform{A}}^{\qquad\! \dform{a}}
     - \Up_{a^1}^{} \X_{A^0\vec{A}}^{\quad \vec{a}} \\
  \wh{A}_{ab^1b^2} 
  =& A_{ab^1b^2} + \Up^p C_{pab^1b^2} \\
  \wh{\na}_{a^1} (C \blacklozenge \si)_{c\dform{a}} 
  =& \na_{a^1}  (C \blacklozenge \si)_{c\dform{a}} 
     + (k-2) \Up_{a^1}  (C \blacklozenge \si)_{c\dform{a}} \\
   & + \bg_{ca^1} \Up^r (C \blacklozenge \si)_{r\dform{a}}
\end{split}
\end{align}
The first two transformations are immediate from
(\ref{transformation}) since $F_{A^0\form{A}}$ is (assumed to be)
conformally invariant.  The next two formulae are directly the
properties of $\Z$-- and $\X$--tractors from (\ref{transformation}).
The last but one 
%follows from the formula (\ref{tractcurv}) for the tractor curvature 
is a simple calculation using the conformal transformation formulae
from for example \cite{GoPetLap},
and the last follows from Lemma \ref{l.lozenge} (i) and
(\ref{grad_trans_(j,k)}). Applying (\ref{trans_proof}) to the formula
(\ref{F_A^0A}) for $\Ph_c$, we obtain
\begin{align*}
  \wh{\Ph}_c & (F_{A^0\form{A}}) - \Ph_c (F_{A^0\form{A}}) 
     = \X_{A^0\form{A}}^{\quad \form{a}} \Bigl[
     - \frac{k+1}{2} \Up^{a^0} C_{a^0a^1c}^{\qquad p} 
       \si_{p\dform{a}}^{} \\
    &- \frac{k(k-1)}{2(n-k)} \Up_{a^1}  
       (C \blacklozenge \si)_{c\dform{a}}
     + \Up^q C_{q\ ca^1}^{\ p} \si_{p\dform{a}}^{}
     - \frac{1}{2} \Up^q C_{q\ a^1a^2}^{\ p} \si_{pc\ddform{a}}^{} \\
    &- \frac{k}{2} C_{a^1a^2c}^{\qquad p} \Up^q \si_{qp\ddform{a}}^{}
     + \frac{k(k\!-\!2)}{2(n\!-\!k)} \Up_{a^1}
       (C \blacklozenge \si)_{c\dform{a}}
     + \frac{k}{2(n\!-\!k)} \bg_{ca^1} \Up^r 
       (C \blacklozenge \si)_{r\dform{a}}
     \Bigr]  
\end{align*}
It is straightforward to verify that sum of the three terms involving
$C \blacklozenge \si$ is equal to 
\begin{equation} \label{prooflozenge}
  - \frac{k}{n-k}\Up^r \bg_{a^1[r}(C \blacklozenge \si)_{c]\dform{a}}~.
\end{equation}
Summing the remaining terms on the right hand side yields
\begin{equation} \label{proofnonlozenge}
\begin{split}
  & \bigl( -\Up^q C_{qa^1c}^{\quad\ p} \si_{p\dform{a}}^{}
      +\frac{k-1}{2} \Up^q C_{a^2a^1c}^{\qquad p} \si_{pq\ddform{a}}^{} 
    \bigr) \\
  & + \Up^q C_{ca^1q}^{\quad\ p} \si_{p\dform{a}}^{}
    - \frac{1}{2} \Up^q C_{a^1a^2q}^{\qquad p} \si_{pc\ddform{a}}^{}
    + \frac{k}{2} \Up^q C_{a^1a^2c}^{\qquad p} \si_{pq\ddform{a}}^{} \\
  & = - \Up^r \left[
        C_{rca^1}^{\quad\ p} \si_{p\dform{a}}^{}
        + C_{a^1a^2[r}^{\qquad p} \si_{|p|c]\ddform{a}}^{}
      \right] ~.
\end{split}
\end{equation}
Now summing the last two displays and comparing the result with 
the definition of $C \lozenge \si$ in (\ref{lozenge}), 
the Lemma (ii) follows.
\end{proof}

We have shown that, in contrast to $\Om_{pcA^0A^1} \si^p$,
$\Phi_c$ for $k \geq 2$ is not conformally invariant. Also note
that it is not algebraic but is rather a first order differential
operator. We would like to replace $\Phi_c$ with an operator which, in a 
suitable sense, has the same essential properties (including linearity) 
and yet is conformally invariant and algebraic.
We deal with invariance first.
For $k \geq 2$, we define the 1st order differential operator
$$ \Ps_c: \cE_{[A^0\form{A}^k]} \longrightarrow \cE_{c[A^0\form{A}^k]}, $$
for a given choice $g \in [g]$ of the metric and 
a section $F_{A^0\form{A}} \in \cE_{[A^0\form{A}^k]}$ (taken to be of the form
(\ref{F_A^0A})), by
\begin{equation} \label{Ps}
  \Ps_c(F_{A^0\form{A}}) := \Ph_c(F_{A^0\form{A}})
  + \frac{1}{n-2} \X_{A^0\form{A}}^{\quad \form{a}} 
    \na^p (C \lozenge \si)_{pc\form{a}}.
\end{equation}
Recall that $(C \lozenge \si)_{[pq]\form{a}} \in \cE(2,k)_0[k+1]$ and
is by construction conformally invariant.  Hence we have the
conformal transformation
$$ \widehat{\na}^p (C \lozenge \si)_{pc\form{a}} 
   = \na^p (C \lozenge \si)_{pc\form{a}} 
   + (n-2) \Up^p (C \lozenge \si)_{pc\form{a}} $$
according to (\ref{grad_trans_(j,k)}). From this and the previous
Lemma (ii) it follows that $\Psi_c$ is conformally invariant.

Now recall we have proved in Lemma \ref{l.whitelozenge} that
$C \lozenge \si =0$ for $\si$ satisfying (\ref{ke}).
Therefore $\Ph_c = \Ps_c$ in this case and we have 

\begin{lem} \label{c.nonclosed+inv}
Lemma \ref{l.nonclosed_noninv} part {\rm (i)} holds if we replace
the operator $\Ph_c$ by $\Ps_c$ therein. \qed
\end{lem} 

Now we replace the operator $\Ps_c$ with an algebraic alternative in
the following way. From (\ref{Ps}) and the formulae (\ref{Ph}) for
$\Phi_c$, it is clear that in the operator $\Ps_c$, 
applied to $F_{A^0\form{A}}$ in the form (\ref{F_A^0A}), 
only the coefficient
of $\X$ contains terms of the first order. Recall that we have the
decomposition $\mathcal{E}_{c\vec{a}^k}[k+1] \cong
\mathcal{E}_{[c\vec{a}^k]}[k+1] \oplus
\mathcal{E}_{\{c\vec{a}^k\}_0}[k+1] \oplus
\mathcal{E}_{\vec{a}^{k-1}}[k-1]$. 
If $\si_{\form{a}} = (k+1)\X_{\quad \form{a}}^{A^0\form{A}}
F_{A^0\form{A}}^{}$ is a solution of (\ref{ke}),
the parts of $\na_c\si_{\form{a}}$
that lie in $\mathcal{E}_{[c\vec{a}^k]}[k+1]$ and
$\mathcal{E}_{\vec{a}^{k-1}}[k-1]$ may be replaced by, respectively,
$\mu_{a^0\form{a}} \in \ce_{a^0\form{a}^k}[k+1]$ and $\nu_{\dform{a}}
\in \ce_{\dform{a}^k}[k-1]$, according to Proposition
\ref{p.prolong_start}. 
Moreover, it is clear that in fact this replacement is 
conformally invariant for any $F_{A^0\form{A}}^{}$.
Thus if we remove, from the $\X$--slot of the
formulae for $\Ps_c$, all the terms depending on $\na_{\{c}
\si_{\form{a}\}_0}$, then the resulting operator $\wt{\Ps}_c$ will be
algebraic, conformally invariant and will satisfy Lemma
\ref{c.nonclosed+inv} (or rather the alternative version of this with
$\wt{\Ps}_c$ replacing $\Ps_c$).
We describe $\wt{\Ps}_c$ explicitly in the following Proposition.

\begin{prop} \label{l.wtPs}
The mapping
$$ \si_{\form{a}} \mapsto {\D}(\si)_{A^0\form{A}}, 
   \quad \mbox{with inverse} \quad
   F_{A^0\form{A}} \mapsto (k+1)\X_{\quad \form{a}}^{A^0\form{A}} 
     F_{A^0\form{A}}^{}~, $$
gives a conformally invariant bijective mapping
between sections of $\si_{\form{a}} \in \cE_{\form{A}^k}[k+1]$
satisfying (\ref{ke}) and sections $F_{A^0\form{A}} \in
\cE_{A^0\form{A}^k}$ satisfying, 
$$\na_c F_{A^0\form{A}} = \wt{\Ps}_c(F_{A^0\form{A}}) \quad \quad 1\leq k\leq n-1~.
$$

For choice $g \in [g]$ of a metric from the conformal class and 
a section $F_{A^0\form{A}} \in \cE_{A^0\form{A}^k}$, expressed  in the form
(\ref{F_A^0A}), the conformally invariant algebraic operator
$\wt{\Ps}_c: \cE_{A^0\form{A}^k} \to \cE_{cA^0\form{A}^k}$
is given by the formula 
\begin{equation} \label{wtPs}
\begin{split}
  \wt{\Ps}_c(F_{A^0\form{A}}) =
  & - \frac{1}{2} \Z_{A^0\form{A}}^{\,a^0\,\form{a}}
      C_{a^0a^1c}^{\qquad\! p} \si_{p\dform{a}}^{}
    + \frac{k(k-1)}{2(n-k)} \W_{A^0A^1\dform{A}}^{\qquad\! \dform{a}}
      (C \blacklozenge \si)_{c\dform{a}}^{} \\
  & + \X_{A^0\form{A}}^{\quad \form{a}} \Bigl[
      A_{a^1c}^{\quad p} \si_{p\dform{a}}^{}
      + \frac{k-1}{2(n-k)} T(\si)_{c\form{a}} \Bigr]
\end{split}
\end{equation}
where 
\begin{align*}
  T(\si)_{c\form{a}}
  =& \frac{1}{2} \bigl( \na_c^{} C_{a^1a^2}^{\quad\ pq} \bigr)
                 \si_{pq\ddform{a}}^{}
     + 2 A^p_{\ ca^1} \si_{p\dform{a}}^{}
     - A^p_{\ a^1a^2} \si_{pc\ddform{a}}^{}
     - \bg_{ca^1} A_{a^2}^{\;\ pq} \si_{pq\ddform{a}}^{}  \\
   & - \bigl( C_{ca^1}^{\quad pq} \mu_{pq\ddform{a}}^{}
       + C_{a^2a^1}^{\quad\ pq} \mu_{pqc\dddform{a}}^{} \bigr) 
     - \frac{n-k-1}{k} C_{a^1a^2c}^{\qquad\! p} \nu_{p\ddform{a}}^{} \\
 \in & \; \cE(1,k)[k-1].
\end{align*}
\end{prop}

\begin{proof} 
The case $k=1$ is just reformulation of Lemma \ref{l.nonclosed_noninv}.
Given Lemma \ref{c.nonclosed+inv}, for the cases $k\geq 2$ 
this boils down to simply 
checking the formula for $\wt{\Ps}$.
This is a direct computation of the formula (\ref{Ps})
for $\Ps_c$
and then in this formula, formally replacing each instance of 
$\na_c \si_{\form{a}}$ by
$ \mu_{c\form{a}} + 
\bg_{ca^1} \nu_{\dform{a}}$. 
We need to compute only the non-algebraic 
terms $\na_{a^1} ( C \blacklozenge \si)_{c\dform{a}}$ from (\ref{Ph})
and $\na^q (C \lozenge \si)_{qc\form{a}}$ from (\ref{Ps}).
The latter is the subject of  Lemma \ref{l.grad_white} below, while the 
the former is dealt with during the proof of that same Lemma, see 
(\ref{grad_black_down}).
Combining these results with (\ref{Ph}) and collecting terms yields the 
formula \nn{wtPs}.
\end{proof}

It remains then to calculate $\na^q (C \lozenge \si)_{qc\form{a}}$ as
required in the proof of the Proposition above.  For this we will need
the following identities. They follow from the (second) Bianchi
identity $\na_{[a}R_{bc]de}=0$ after a short computation.
\begin{equation} \label{Bianchi2}
\begin{split}
  \na_{a^1} C_{ca^2b^1b^2}
  &= \frac{1}{2} \na_c C_{a^1a^2b^1b^2} 
     - \bg_{cb^1}A_{b^2a^1a^2} 
     + 2 \bg_{a^1b^1} A_{b^2ca^2} \\
  \na_{a^1} C_{a^2a^3b^1b^2} 
  &= 2 \bg_{a^1b^1} A_{b^2a^2a^3}.
\end{split}
\end{equation}

\begin{lem} \label{l.grad_white}
Assume $ 2\leq k\leq n-1$.
If the $\si_{\form{a}} \in \cE_{\form{a}^k}[k+1]$ then, up to the addition of (conformally invariant) terms involving the Weyl curvature contracted into $\na_{\{c}
\si_{\form{a}\}_0}$,
$\na^q (C \lozenge \si)_{qc\form{a}} \in \cE(1,k)_0[k-1]$ 
is given by the formula
\begin{align*}
& \frac{n-2}{2(n-k)} \Bigl[
    \frac{1}{2} (\na_c^{} C_{a^1a^2}^{\quad\ pq}) \si_{pq\ddform{a}}^{}
    - \bigl( C_{ca^1}^{\quad pq} \si_{pq\ddform{a}}^{}
      \!+\! C_{a^2a^1}^{\quad\ pq} \si_{pqc\dddform{a}}^{} \bigr) \\
  & \! + (n-k-1) \bigl( A^p_{\ a^1a^2} \si_{pc\ddform{a}}^{}
           + 2 A^p_{\ a^1c} \si_{p\dform{a}}^{} \bigr)
    + \frac{(n-k+1)}{k} C_{a^1a^2c}^{\qquad p} \nu_{p\ddform{a}}^{} \\
  & \! + \frac{(k-2)}{k} \bg_{ca^1}^{} C_{a^2a^3}^{\quad\ pq}
         \nu_{pq\dddform{a}}^{}
    - (k-1) \bg_{ca^1}^{} A_{a^2}^{\ \ pq} \si_{pq\ddform{a}}^{} \Bigr] 
    + (n\!-\!2) A_{a^1c}^{\quad p} \si_{p\dform{a}}^{}. 
\end{align*}
\end{lem}

\begin{proof}
Here we simply expand $\na^q (C \lozenge \si)_{qc\form{a}} $ via the
Leibniz rule and in the process we will formally replace each $\na_c
\si_{\form{a}}$ by $\mu_{c\form{a}} + \bg_{ca^1} \nu_{\dform{a}}$.
We shall start with
$\na_{a^1} (C \blacklozenge \si)_{c\dform{a}}$.
Recall $(C \blacklozenge \si)_{c\dform{a}}$ was given in 
(\ref{lozenge}) as a sum of two terms. 
Applying $\na_{a^1}$ to these, we obtain
\begin{align*}
  \na_{a^1}^{} C_{ca^2}^{\quad pq} \si_{pq\ddform{a}}^{}
  =& \frac{1}{2} (\na_c^{} C_{a^1a^2}^{\quad\ pq}) \si_{pq\ddform{a}}^{}
     - A^q_{\ a^1a^2} \si_{cq\ddform{a}}^{}
     + 2 A^q_{\ ca^2} \si_{a^1q\ddform{a}}^{} \\
   & + C_{ca^2}^{\quad pq} \bigl( \mu_{a^1pq\ddform{a}} 
       + \bg_{a^1[p} \nu_{q\ddform{a}]}^{}) \\
  \na_{a^1}^{} C_{a^3a^2}^{\quad\ pq} \si_{pqc\dddform{a}}^{}
  =& 2 A^q_{\ a^3a^2} \si_{a^1qc\dddform{a}}^{} 
     + C_{a^3a^2}^{\quad\ pq} \bigl( \mu_{a^1pqc\dddform{a}}
       +  \bg_{a^1[p} \nu_{qc\dddform{a}]}^{}). 
\end{align*}
where we have also used (\ref{Bianchi2}).
Now summing of the right-hand sides of the last displays yields
\begin{equation} \label{grad_black_down}
\begin{split}
  & \na_{a^1} ( C \blacklozenge \si)_{c\dform{a}}
  = \frac{k\!-\!2}{k} \Bigl[  
    \frac{1}{2} (\na_c^{} C_{a^1a^2}^{\quad\ pq}) \si_{pq\ddform{a}}^{}
    - A^p_{\ a^1a^2} \si_{pc\ddform{a}}^{}
    + 2 A^p_{\ ca^1} \si_{p\dform{a}}^{} \\
  & -\! \bigl( C_{ca^1}^{\quad pq}  \mu_{pq\dform{a}}^{}
      \!+\! C_{a^2a^1}^{\quad\ pq}  \mu_{pqc\ddform{a}}^{}  \bigr)
    + \frac{1}{k} C_{a^1a^2c}^{\quad\ \ p} \nu_{p\ddform{a}}^{} 
    - \frac{1}{k} \bg_{ca^1} C_{a^2a^3}^{\quad\ pq} \nu_{pq\dddform{a}}^{}
    \Bigr]
\end{split}
\end{equation}
where we have used $\frac{2}{k} C_{ca^2a^1}^{\quad\ \ q} = 
\frac{1}{k} C_{a^1a^2c}^{\quad\ \ q}$. Note 
$\na_{a^1} ( C \blacklozenge \si)_{c\dform{a}} \in \cE(1,k)[k-1]$.

Now we shall compute the formula for $\nd^q (C \lozenge \si)_{qc\form{a}}$.
According to (\ref{lozenge}), $(C \lozenge \si) $ is defined
as sum of three terms. Applying $\na^q$ to the first of these,
and using \nn{bi1}, we obtain 
$$ \na^q C_{qca^1}^{\quad p} \si_{p\dform{a}}^{} 
   = (n-3) A_{ca^1}^{\quad p} \si_{p\dform{a}}^{}
     + C_{\ ca^1}^{q \ \ p} (\mu_{qp\dform{a}}^{} 
                             + \bg_{q[p} \nu_{\dform{a}]}^{})~. $$
Similarly for the second term, we obtain
\begin{align*}
  \na^q C_{a^1a^2[q}^{\qquad p} \si_{|p|c]\ddform{a}}^{} 
  =& \frac{1}{2} (n-3) A^p_{\ a^1a^2} \si_{pc\ddform{a}}^{}
     + \frac{1}{2} C_{a^1a^2}^{\quad\ qp} 
       ( \mu_{qpc\ddform{a}}^{} + \bg_{q[p}\nu_{c\ddform{a}]}^{} ) \\
   & - \frac{1}{2} (\na^q C^{\ p}_{c\ a^1a^2}) \si_{pq\ddform{a}}^{}
     + \frac{n-k+1}{2k} C_{a^1a^2c}^{\qquad p} \nu_{p\ddform{a}}^{} ~,
\end{align*}
where we have used $\na^q\si_{q\dform{a}} = \frac{n-k+1}{k} \nu_{\dform{a}}$.
Summing the right hand sides of the last two displays with 
the third term $\frac{k}{n-k}\na^q \bg_{a^1[q} 
(C \blacklozenge \si)_{c]\dform{a}}$ 
yields
\begin{equation} \label{formula1}
\begin{split}
  \na^q (C & \lozenge \si)_{qc\form{a}} 
    = \frac{1}{2} (\na^p C^{\ q}_{c\ a^1a^2}) \si_{pq\ddform{a}}^{}   
    - \frac{1}{2} \bigl( C_{ca^1}^{\quad pq} \mu_{pq\dform{a}}^{}
      \!+\! C_{a^2a^1}^{\quad\ pq} \mu_{pqc\ddform{a}}^{}  \bigr) \\
%    - \frac{k+1}{2(k-1)} (C \lozenge \mu)_{c\form{a}} \\
  & + (n-3) \Bigr[ A_{ca^1}^{\quad p} \si_{p\dform{a}}^{}
      + \frac{1}{2} A^p_{\ a^1a^2} \si_{pc\ddform{a}}^{} \Bigr]
    + \frac{n-1}{2k} C_{a^1a^2c}^{\qquad p} \nu_{p\ddform{a}}^{} \\
  & + \frac{k}{2(n-k)} \na_{a^1} (C \blacklozenge \si)_{c\dform{a}}
    - \frac{k}{2(n-k)} \bg_{ca^1} \na^q (C \blacklozenge \si)_{q\dform{a}}
\end{split}
\end{equation}
where we have used $C_{\;\ ca^1}^{[q \ \ p]} = - \frac{1}{2}
C_{ca^1}^{\quad qp}$. In the last display, we need the term
$\na^p (C \blacklozenge \si)_{p\dform{a}}$. Using the definition
(\ref{lozenge}) and applying the Leibniz rule for $\na^p$,
we obtain
\begin{equation} \label{grad_black_up}
\begin{split}
  \na^p (C \blacklozenge \si )_{p\dform{a}}
  =& \frac{k-2}{k} \Bigl[
     (n-3) A_{a^2}^{\ \ pq} \si_{pq\ddform{a}}^{}
     + C^{r\;\ pq}_{\ a^2} \bg_{r[p}^{} \nu_{q\ddform{a}]}^{} \\
   & + (\na^r_{} C^{pq}_{\quad a^3a^2}) \si_{pqr\dddform{a}}
     - \frac{n-k+1}{k} C_{a^2a^3}^{\quad\ pq} \nu_{pq\dddform{a}}
     \Bigr] \\
  =& \frac{(k-2)(n-1)}{k} \Bigl[
     A_{a^2}^{\ \ pq} \si_{pq\ddform{a}}^{}
     - \frac{1}{k} C_{a^2a^3}^{\quad\ pq} \nu_{pq\dddform{a}} \Bigr]
\end{split}
\end{equation}
using (\ref{Bianchi2}). We will also need the identity
$$ \frac{1}{2} (\na^p C^{\ q}_{c\ a^1a^2}) \si_{pq\ddform{a}}^{}
   = + \frac{1}{4} (\na_c^{} C^{pq}_{\quad a^1a^2}) \si_{pq\ddform{a}}^{}
     - \frac{1}{2} \bg_{ca^1} A_{a^2}^{\;\ pq} \si_{pq\ddform{a}}^{}
     + A_{a^1c}^{\quad p} \si_{p\dform{a}}^{} $$
which uses (\ref{Bianchi2}). Now we are ready to simplify (\ref{formula1})
using (\ref{grad_black_down}), (\ref{grad_black_up}) and the last display.
Collecting terms 
the result is
\begin{align*}
  &\na^q (C \lozenge \si)_{qc\form{a}}
  = \frac{n-2}{4(n-k)} \Bigl[    
    (\na_c^{} C_{a^1a^2}^{\quad\ pq}) \si_{pq\ddform{a}}^{}
    - 2 \bigl( C_{ca^1}^{\quad pq} \mu_{pq\dform{a}}^{}
      \!+\! C_{a^2a^1}^{\quad\ pq} \mu_{pqc\ddform{a}}^{} \bigr) \\
  & \qquad\qquad\qquad 
    + 2(n-k-1) A^p_{\ a^1a^2} \si_{pc\ddform{a}}^{}
    + \frac{2(n-k+1)}{k} C_{a^1a^2c}^{\qquad p} \nu_{p\ddform{a}}^{} \\
  & \qquad\qquad\qquad 
    + \frac{2(k-2)}{k} \bg_{ca^1}^{} C_{a^2a^3}^{\quad\ pq} 
      \nu_{pq\dddform{a}}^{}
    - 2(k-1) \bg_{ca^1}^{} A_{a^2}^{\ \ pq} \si_{pq\ddform{a}}^{} \Bigr] \\
  & + \frac{1}{(n-k)} \Bigl[ (n-k) A_{a^1c}^{\quad p} 
    + (k-2) A^p_{\ ca^1} 
    + (n-3)(n-k) A_{ca^1}^{\quad p}
    \Bigr] \si_{p\dform{a}}^{} 
\end{align*}
Now the final step is to simplify the last line using the relation
$A_{ca^1}^{\quad p} = A^p_{\ a^1c} + A_{a^1c}^{\quad p}$ which follows
directly from the definition  $A_{pa^1c} := 2 \na_{[a^1} P_{c]p}$. A short computation
reveals that the last line is equal to
$$ (n-2) A_{a^1c}^{\quad\, p} + (n-2) \frac{n-k-1}{n-k} A^p_{\ a^1c}. $$
The Lemma now follows from the last two displays.
\end{proof}

Summarising our results we have the following.
\begin{thm} \label{main}
For $1\leq k\leq n-1$, the mapping $\cE_{\form{a}^k}[k+1] \to
\cE_{A^0\form{A}^k}$ given by $\si \mapsto \D(\si)$ defined by
(\ref{si_A}) is a conformally invariant differential operator.  Upon
restriction it gives a bijective mapping from solutions of the
conformal Killing equation \nn{ke} onto sections of
$\cE_{A^0\form{A}^k}$ that are parallel with respect to the connection
${}^k\na_c := \na_c - \wt{\Ps}_c$ where $\na_c$ is the normal tractor
connection and $\wt{\Ps}_c$ is given by (\ref{wtPs}). The connection
${}^k\na_c$ is a conformally invariant connection on the form-tractor
bundle $\cE_{A^0\form{A}^k}$. The inverting map from sections of
$\cE_{A^0\form{A}^k}$, parallel for ${}^k\na_c$, to solutions of
\nn{ke} is $F_{A^0\form{A}} \mapsto (k+1)\X_{\quad
\form{a}}^{A^0\form{A}} F_{A^0\form{A}}^{}$.

Sections of
$\cE_{A^0\form{A}^k}$ which are parallel for the normal tractor
connection $\nd_c$ are mapped injectively to solutions of \nn{ke} by
$$
F_{A^0\form{A}} \mapsto (k+1)\X_{\quad \form{a}}^{A^0\form{A}} 
     F_{A^0\form{A}}^{}~,  
$$
and $\tilde\Psi_c$ annihilates the range of this map. 
\end{thm}
\begin{proof} Everything has been established in the previous Lemmas 
except for the last claim. That parallel sections are mapped
injectively to conformal Killing forms is an immediate consequence of
the formula \nn{normtrconn} for the normal tractor connection on
form-tractors.  (Note that the equation from the first slot of $\na_c
F_{A^0\form{A}} =0$ is $\na_c \si_{\vec{a}^k} - (k+1) \mu_{c\vec{a}^k}
+ \bg_{ca^1} \ph_{\dot{\vec{a}}^k}=0$. This is the same equation as
from the first slot for a $(k+1)$-form-tractor parallel for
${}^k\na_c$, as $\wt{\Ps}_c$ does not affect this top slot -- the
coefficient of $\Y$.) 
Next it is an
elementary exercise using the formula \nn{normtrconn} 
to verify that if
$F_{A^0\form{A}}$ is parallel for the normal tractor connection, then 
necessarily $F_{A^0\form{A}} = \D(\si)$ where $\si_{\form{a}}=(k+1)\X_{\quad
\form{a}}^{A^0\form{A}} F_{A^0\form{A}}^{}$.
On the other hand from the first part of the
Theorem it follows that $\D(\si)$ is parallel for ${}^k\na$. So 
$\tilde\Psi_c(\si)$ vanishes everywhere. 
\end{proof}
\noindent {\bf Remark:} Let us say (as suggested in \cite{leit}) that 
a conformal Killing form $\si$ is {\em normal} if it has the property
that $\D(\si)$ is parallel for the normal tractor connection. It
follows immediately from the Theorem that the operator $\tilde\Psi_c$
detects exactly the failure of conformal Killing forms to be normal; a
conformal Killing form is normal if and only if $\tilde\Psi_c(\si)$ is
zero.

\vspace{1ex}

If $\si\in \ce^k[k+1]$ vanishes on an open set then note that $\D
(\si)$ vanishes on the same open set since $\D $ factors through the
universal jet operator $j^2$. On the other hand if $\si$ is a
conformal Killing form then, from the Theorem $\D(\si)$ is parallel
for the connection ${}^k\na$. Thus we have the following.
\begin{cor}
On connected manifolds $M$ a non-trivial conformal Killing
form is non-vanishing on an open dense subspace.
\end{cor}
\noindent The corollary here does not use the conformal invariance of
the connection and so this conclusion also follows from
\cite{Sem-MathZ}.

\section{Coupled conformal Killing equations} \label{s.helicity}

  In this section we show that solutions $\si\in \ce^k[k+1]$ of the
  original equation (\ref{ke}) are in bijective correspondence with
  solutions of the coupled conformal Killing equation
  $\tilde{\nd}_{(a}\ol{\si}_{b)_0\vec{B}^{k-1} }=0$ on
  $\mathcal{E}_{a\vec{B}^{k-1}}[2]$ for a certain conformally
  invariant connection $\tilde{\na}$.   Along the way we obtain
  some related preliminary results that should be of independent
  interest.

First let us observe that for any form $\si\in \ce^k[k+1]$, $1\leq
k\leq n-1$, we may form the tractor-valued forms
\begin{equation}\label{ulol}
 \overline{\si}_{\vec{a}^{k-l}\vec{B}^l} =
\ol{M}_{\vec{B}^l}^{\vec{a}^{k-l,l}} \si_{\vec{a}^k} \quad \mbox{and} \quad
\underline{\si}_{\vec{a}^{k+l}\vec{B}^l} = \ul{M}_{\vec{a}^{k,l}\vec{B}^l}
\si_{\vec{a}^k}
\end{equation}
 where the invariant differential splitting operators
$\ol{M}$ and $\ul{M}$ are defined by the formulae, for $1\leq l \leq k$, 
\begin{gather*}
  \ol{M}_{\vec{B}^l}^{\vec{a}^{k-l,l}}
    : \mathcal{E}_{\vec{a}^k}[k+1] \longrightarrow 
       \mathcal{E}_{\vec{a}^{k-l}\vec{B}^l}[k-l+1] \\
  \ol{M}_{\vec{B}^l}^{\vec{a}^{k-l,l}} \si_{\vec{a}^k}
    = (n-k+1) \Z_{\vec{B}^l}^{\vec{b}^l} \si_{\vec{a}^{k-l}\vec{b}^l}
       - l    \X_{B^1\vec{\dot{B}}^l}^{\quad \vec{\dot{b}}^l} \na^{b^1} 
              \si_{\vec{a}^{k-l}\vec{b}^l} 
\end{gather*}
and, for $1\leq l\leq n-k$, 
\begin{gather*}
  \ul{M}_{\vec{a}^{k,l}\vec{B}^l} : \mathcal{E}_{\vec{a}^k}[k+1]
    \longrightarrow \mathcal{E}_{\vec{a}^{k+l}\vec{B}^l}[k+l+1] \\
    \ul{M}_{\vec{a}^{k,l}\vec{B}^l} \si_{\vec{a}^k} = (k+1)
    \Z_{\vec{B}^l}^{\vec{b}^l} \bg_{\vec{b}^l\vec{a}^{k,l}}
    \si_{\vec{a}^k} - l \X_{B^1\vec{\dot{B}}^l}^{\quad
    \vec{\dot{b}}^l} \bg_{\vec{\dot{b}}^l\vec{\dot{a}}^{k,l}}
    \na_{a^{k+1}} \si_{\vec{a}^k} ~.\\
\end{gather*}
Here we use multi-indices  
$$
\begin{aligned}
&\vec{a}^{k,l} =[a^{k+1} \cdots a^{k+l}]\\
&\vec{\dot{a}}^{k,l}= [a^{k+2} \cdots a^{k+l}]~.
\end{aligned}
$$ The conformal invariance of $\ol{M}$ and $\ul{M}$ may be verified
directly via the formulae (\ref{transformation}).  

Although $\ol{\si}_{\vec{a}^{k-l}\vec{B}^l}$ and
$\ul{\si}_{\vec{a}^{k+l}\vec{B}^l}$, as defined in \nn{ulol}, are
invariant for the stated ranges of $l$, in the sequel we shall only need 
the tensor valence of $\ol{\si}$ and $\ul{\si}$ to be in the interval
$[1,n-1]$. Therefore we shall henceforth assume that for
$\ol{\si}_{\vec{a}^{k-l}\vec{B}^l}$ we have $1\leq l \leq k-1$ and for
$\ul{\si}_{\vec{a}^{k+l}\vec{B}^l}$ we have $1\leq l\leq n-k-1$,
respectively. 

Let us next describe 
$\na_{\{c} \ol{\si}_{\vec{a}^{k-l}\}_0\vec{B^{l}}}$ and $\na_{\{c}
\ul{\si}_{\vec{a}^{k+l}\}_0\vec{B^{l}}}$ when $\si$ is a solution
of (\ref{ke}). (Recall that $\na$ denotes the coupled
Levi--Civita--normal tractor connection.)  This is explicitly
formulated in the proposition below. First we need the following lemma.

\begin{lem} \label{l.helicity}
Let us suppose that $\si$ is a solution of (\ref{ke}). Then
\begin{align}
  \na_c\na^p \si_{\vec{a}^{k-l}p\dot{\vec{b}}^l}
    &\stackrel{{\{c\vec{a}^{k-l}\}_0}}{=} (n-k+1) \left[
      -\frac{k-1}{n-k} C_{c\ [a^1}^{\ p\;\ q} 
      \si_{|p|\dot{\vec{a}}^{k-l}|q|\dot{\vec{b}}^l]}
      -P_c^{\ p} \si_{\dot{\vec{a}}^{k-l}p\dot{\vec{b}}^l} \right] 
    \tag{a} \\
  \na_c\na_{a^{k+1}}\si_{\vec{a}^k}
    &\stackrel{{\{c\vec{a}^{k+1}\}_0}}{=} (k+1) \left[
       C_{ca^{k+1}a^1}^{\quad\qquad p} \si_{p\dot{\vec{a}}^k} 
       - P_{ca^{k+1}} \si_{\vec{a}^k}
    \right]. 
    \tag{b}
\end{align}
\end{lem}
\noindent In reading (b) here recall the convention that sequentially 
labelled indices (at a given level) are assumed to be skewed over. 

\begin{proof} 
First let us note that the trace part in the first case, and
skew--symmetrisation $[c\vec{a}^{k+1}]$ in the second case, is zero on
both sides. In the subsequent discussion we use Proposition
\ref{p.prolong_start} and the notation therein.

The left-hand side of (a) is equal to $\frac{n-k+1}{k} \nd_c
\nu_{\form{a}^{k-l}\dform{b}^l}$ up to
the sign $(-1)^{k-l}$.  Now the Lemma (a) follows using $C_{c\;\
a^1}^{\ [p\;\ q]} = \frac{1}{2} C_{ca^1}^{\quad pq}$ and the equation
for $\na_c^{} \nu_{\form{a}^{k-l}\dform{b}^l}$ in
(\ref{prolong_start}) where $(C \blacklozenge
\si)_{c\form{a}^{k-l}\dform{b}^l}$ is given by Lemma \ref{l.lozenge}
(i). Note that the projection $\{..\}$ over indices in the latter
lemma exactly removes the completely skew--symmetric part of
$C_{ca^2}^{\quad pq} \si_{pq\ddform{a}}^{}$ (see \nn{black}). Since the
projection $\{c\vec{a}^{k-l}\}_0$ annihilates the completely
skew--symmetric part $C_{[ca^2}^{\quad\ pq} \si_{|pq|\ddform{a}]}^{}$
we have 
$(C \blacklozenge \si)_{c\form{a}^{k-l}\dform{b}^l}
=_{{\{c\vec{a}^{k-l}\}_0}} 
C_{ca^1}^{\quad pq} \si_{pq \dform{a}^{k-l}\dform{b}}$.  
The part (b) follows similarly from the
expression for $\na_c \mu_{a^{k+1}\form{a}^k}$ in
(\ref{prolong_start}).
\end{proof}

\begin{prop} \label{charac}
The form $\si \in \cE^k[k+1]$, $1 \leq k \leq n-1$ 
is a solution of 
(\ref{ke}) if and only if either of the following
conditions is satisfied: 
\begin{align*}
\na_c \ol{\si}_{\vec{a}^{k-l}\vec{B}^l}
  & \stackrel{{\{c\vec{a}^{k-l}\}_0}}{=} 
\frac{l(k-1)(n-k+1)}{n-k} \X_{B^1\dot{\vec{B}}^l}^{\quad \dot{\vec{b}}^l}
     C_{c\ [a^1}^{\ p\;\ q} 
     \si_{|p|\dot{\vec{a}}^{k-l}|q|\dot{\vec{b}}^l]} \\
\na_c \ul{\si}_{\vec{a}^{k+l}\vec{B}^l}
  &\stackrel{{\{c\vec{a}^{k+l}\}_0}}{=} -l(k+1)  
     \X_{B^1\dot{\vec{B}}^l}^{\quad \dot{\vec{b}}^l}
     C_{c[a^{k+1}a^1}^{\quad\qquad p} 
     \si_{|p|\dot{\vec{a}}^k} \bg_{\dot{\vec{a}}^{k,l}]\dot{\vec{b}}^l}.
\end{align*}
\end{prop}

\begin{proof}  
The expressions
on the left-hand-side can be computed by directly differentiating the
expressions \nn{ulol} defining $\ul{\si}$ and $\ol{\si}$ and expanding
in terms of the $\X$, $\Y$, $\W$, $\Z$ splitting operators introduced
in \ref{s.form_tractor}. The resulting ``$\Y$--slot'' (i.e. the
coefficient of $\Y$) on the left-hand-side is zero order, as an
operator on $\si$, and is killed by the symmetrisation
$\{c\vec{a}^{k-l}\}$ in the case of $\na_c
\ol{\si}_{\vec{a}^{k-l}\vec{B}^l}$ and by taking the trace-free part
in the case of $\na_c \ul{\si}_{\vec{a}^{k+l}\vec{B}^l}$.  Essentially
the same argument shows (in both cases) that also the operator in the
$\W$ slot vanishes.  The $\Z$ slot is of the first order as an
operator on $\si$. To show this vanishes requires some computation. We
will need the relation
\begin{equation} \label{skew_long}
  k \bg_{c[a^1}^{} \na^p \si_{|p|\dot{\vec{a}}^{k-l}\vec{b}^l]} = 
  (k-l) \bg_{ca^1}^{} \na^p \si_{|p|\dot{\vec{a}}^{k-l}\vec{b}^l}
  + l \bg_{cb^1}^{} \na^p \si_{\vec{a}^{k-l}p\dot{\vec{b}}^l}.
\end{equation}
(Recall our convention that all sequentially labelled indices are
implicitly skewed over. So the $b $-indices are skewed and also the
$a$-indices are skewed.) 
To prove this first observe the projection to
the completely skew part of the right-hand-side obviously yields
exactly the left-hand-side. On the other hand the right-hand-side is
manifestly skew over the $b$--indices and also over the
$a$--indices. A trivial calculation verifies that that it is also
skew--symmetric in the index pair $a^1b^1$ and so the result follows.

Using (\ref{prolong_start}) for $\na_c\si_{\form{a}}$, 
it is straightforward
to compute the $\Z$ slot of $\na_c \ol{\si}_{\vec{a}^{k-l}\vec{B}^l}$ is
$$ % \Z_{B^0\vec{B}^l}^{\,b^0\,\vec{b}^l} \Bigl[
     (n-k+1) \na_{[c} \si_{\vec{a}^{k-l}\vec{b}^l]}
     + k \bg_{c[a^1} \na^p \si_{|p|\dot{\vec{a}}^{k-l}\vec{b}^l]}
     - l \bg_{cb^1}^{} \na^p \si_{\vec{a}^{k-l}p\dot{\vec{b}}^l}. $$
The first term is killed by the projection $\cP_{\{c\vec{a}^{k-l}\}}$ and 
the remaining part is in the trace part over
$\{c\vec{a}^{k-l}\}$ (i.e. in particular is annihilated by $\cP_{\{c\vec{a}^{k-l}\}_0}$ )
due to (\ref{skew_long}). The $\Z$ slot of 
$\na_c \ul{\si}_{\vec{a}^{k+l}\vec{B}^l}$ is
$$ %\Z_{B^0\vec{B}^l}^{\,b^0\,\vec{b}^l} \Bigl[ (k+1)
     \bg_{\vec{b}^l\vec{a}^{k,l}} \na_{[c} \si_{\vec{a}^k]} - l
     \bg_{cb^1} \bg_{\dot{\vec{b}}^l\dot{\vec{a}}^{k,l}} \na_{a^{k+1}}
     \si_{\vec{a}^k} + \frac{k(k\!+\!1)}{n\!-\!k\!+\!1} \bg_{ca^1}
     \bg_{\vec{b}^l\vec{a}^{k,l}} \na^p \si_{p\dot{\vec{a}}^k} $$
     (also using (\ref{prolong_start})). The last term is killed by 
     taking the
     trace--free part and it is easy to show the sum of the first two
     terms is $\bg_{\vec{b}^l\vec{a}^{k,l}} \na_c \si_{\vec{a}^k}$
     (up to a scalar multiple) which vanishes after the symmetrisation
     $\{c\vec{a}^{k+l}\}$.

At this point it is worthwhile noting that if the projection
$\cP_{\{c\vec{a}^{k-l}\}_0}$ kills $\na_c
\ol{\si}_{\vec{a}^{k-l}\vec{B}^l}$ or the projection
$\cP_{\{c\vec{a}^{k+l}\}_0}$ kills $\na_c
\ul{\si}_{\vec{a}^{k+l}\vec{B}^l}$ then $\si$ is a solution of
(\ref{ke}); the vanishing of the $\Z$--slots implies 
$\na_c \si_{\form{a}} = \mu_{c\form{a}} + \bg_{ca^1} \nu_{\dform{a}}$
in (\ref{prolong_start}) 
since $\cP_{\{c\vec{a}^{k}\}_0}\circ\cP_{\{c\vec{a}^{k-l}\}_0}$ is a
non-zero multiple of $\cP_{\{c\vec{a}^{k}\}_0}$.

It remains to evaluate the $\X$--slots. This can be done easily using the rules 
for $\na_c \W$ and $\na_c \X$ from \ref{s.form_tractor}. We get
\begin{align*}
-l& \X_{B^1\dot{\vec{B}}^l}^{\quad \dot{\vec{b}}^l} \left[
    (n-k+1) P_c^{\ p} \si_{\vec{a}^{k-l}p\dot{\vec{b}}^l} +
    \na_c \na^p \si_{\vec{a}^{k-l}p\dot{\vec{b}}^l} \right] \\
-l& \X_{B^1\dot{\vec{B}}^l}^{\quad \dot{\vec{b}}^l} \left[
    (k+1) P_{c[a^{k+1}} \si_{\vec{a}^k} g_{\dot{\vec{a}}^{k,l}]\dot{\vec{b}}^l} +
    \na_c \na_{[a^{k+1}} \si_{\vec{a}^k} g_{\dot{\vec{a}}^{k,l}]\dot{\vec{b}}^l} \right]
\end{align*}
for $\na_c \ol{\si}_{\vec{a}^{k-l}\vec{B}^l}$ and $\na_c
\ul{\si}_{\vec{a}^{k+l}\vec{B}^l}$, respectively. Now the proposition follows
using Lemma \ref{l.helicity}.
\end{proof}

For our next construction we will especially need the first case of
the proposition above for $l=k-1$, that is for 
$\ol{\si}_{a^1\dot{\vec{B}}^k}$.  We will construct a
connection $\tilde{\na}$ on $\mathcal{E}_{a^1\vec{\dot{B}}^k}$ such
that the equation
$\tilde{\na}_{(c}^{}\ol{\si}_{a^1)_0\dot{\vec{B}}^k}=0$ is
equivalent to the equation (\ref{ke}). Reformulating the Proposition
for $\ol{\si}_{a^1\dot{\vec{B}}^k}$, we get that $\si$ is a
solution of (\ref{ke}) if and only if
\begin{equation}\label{star} 
\na_{(c}^{} \ol{\si}_{a^1)_0\vec{\dot{B}}^k} =
   \frac{(k-1)(k-2)(n-k+1)}{n-k} \X_{B^2\vec{\ddot{B}}^k}^{\quad
   \vec{\ddot{b}}^k} C_{b^3\ (c}^{\;\ p\;\
   q}\si_{a^1)_0pq\dddot{b}^k}^{}. 
\end{equation}
This shows that $\na_{(c}^{} \ol{\si}_{a^1)_0\vec{\dot{B}}^k}=0$ is equivalent to (\ref{ke}) 
in the flat case. In the curved case we modify the connection $\na$ in the
following way. Let us consider the tensor-tractor field
\begin{align*}
  \kappa_{cE^0E^1F^0F^1}
    :&= \X_{E^0E^1}^{\quad e^1} \Om_{ce^1F^0F^1} \\
     &= \X_{E^0E^1}^{\quad e^1}\Z_{F^0F^1}^{\,f^0f^1} C_{ce^1f^0f^1}
        - 2 \X_{E^0E^1}^{\quad e^1}\X_{F^0F^1}^{\quad f^1} A_{f^1ce^1}~,
\end{align*}
where $\Om_{ce^1F^0F^1}$ is the curvature of the normal tractor
connection. By construction this is conformally invariant.  We will
show that the required connection $\tilde{\na}$ can be written in the
form
$$ \tilde{\na}_c = \na_c + x \kappa_{c}\sharp\sharp, \ x\InR $$ where
(via the tractor metric) we view $\kappa_{cE^0E^1F^0F^1}$ as a 1-form
taking values in ${\rm End}(\ce^A)\otimes {\rm End}(\ce^A)$ and
$\sharp$ indicates the usual action of tractor-bundle endomorphisms
(i.e.\ it is the tractor bundle analogue of the ${\rm End}(TM)$ action
defined in section \ref{fullpro} and we use the same notation as for
that case).  To determine the parameter $x\InR$, let us compute the
double action:
\begin{align*}
  \kappa_{c}\sharp\sharp (\ol{\si}_{a^1\vec{\dot{B}}^k}) &=
  \X^{e^1}\Z^{f^0f^1} C_{ce^1f^0f^1} \sharp\sharp \left[ (n-k+1)
  \Z_{\vec{\dot{B}}^k}^{\,\vec{\dot{b}}^k}
  \si_{\vec{a^1\vec{\dot{b}}^k}^k}\right] \\ 
&= (k-1)(n-k+1)
  \X^{ e^1}\sharp
  \Z_{\vec{\dot{B}}^k}^{\,\vec{\dot{b}}^k} C_{ce^1b^2}^{\qquad q}
  \si_{a^1q\vec{\ddot{b}}^k}^{} \\ &= -\frac{1}{2} (k-1)(k-2)(n-k+1)
  \X_{B^2\vec{\ddot{B}}^k}^{\quad \vec{\ddot{b}}^k} C_{c\ b^3}^{\ p\
  q} \si_{a^1qp\vec{\dddot{b}}^k}^{}.
\end{align*}
The form of the right-hand-side shows that $\tilde{\na}$ is the
required connection for a suitable parameter $x\InR$, and comparing with 
\nn{star} yields the  explicit value for $x$. The
resulting connection is
\begin{equation} \label{mod_connection}
  \tilde{\na}_c = \na_c + \frac{2}{n-k}\kappa_c\sharp\sharp,
\end{equation}
where on the right-hand side $\nd$ is the usual tractor connection.  
 Note that this connection is obviously conformally invariant (since
 both $\kappa$ and the the tractor connection are conformally
 invariant).  This might seem inevitable, since from its derivation
 (or otherwise) it is clear that the equation \nn{star} is conformally
 invariant. However \nn{mod_connection} is an invariant connection
 which may turn out to have applications in other circumstances.

Let us summarise the last result. 
\begin{prop} \label{coupKill}
A weighted $k$-form $\si\in \ce^k[k+1]$ is a conformal Killing
$k$-form (i.e.\ solution of \nn{ke}) if and  only if
\begin{equation}\label{origtw}
\tilde{\nd}_{(a}\ol{\si}_{b)_0}=0
\end{equation}
where $\tilde{\nd}$ is the Levi-Civita connection coupled with
\nn{mod_connection} and $\ol{\si}$ is the conformally invariant
tractor extension of $\si$ given by \nn{ulol} with $l=k-1$.
\end{prop}

Although we shall not directly need it below it is interesting to
observe at this point that the last result generalises. 
First observe that as well as  the action $\kappa_c\sharp\sharp$
used in (\ref{mod_connection}), we can consider also the action
$\om_c\sharp\sharp$ where we view the tensor-tractor field
$$ \om_{cE^0E^1f^0f^1} := \X_{E^0E^1}^{\quad e^1} C_{ce^1f^0f^1} $$ as
a one form taking values in ${\rm End}(\ce^A)\otimes {\rm End}(\ce^a)$
and $\sharp$ indicates the usual action of tensor/tractor-bundle
endomorphisms.  Now for any real parameter $x$ we obtain a connection
on tensor tractor fields via the formula,
\begin{equation}\label{dcoup}
 {\na}_c^x = \na_c + x ( \om_c\sharp\sharp + \ka_{c}\sharp\sharp). 
\end{equation}
where $\nd$ indicates the usual coupled tractor-Levi Civita connection.

\begin{thm} \label{coupKillthm}
A weighted $k$-form $\si\in \ce^k[k+1]$ is a conformal Killing
$k$-form (i.e.\ solution of \nn{ke}) if and  only if either of the 
following conditions holds:
$$
{\na}_{\{c}^x \ol{\si}_{\vec{a}^{k-l}\}_0\vec{B}^l}^{} =0 
\quad \mbox{or} \quad
{\na}_{\{c}^y \ul{\si}_{\vec{a}^{k+l}\}_0\vec{B}^l}^{}=0
$$ 
where $x=\frac{2}{n-k}$ and $y=\frac{2}{k}$,  and $\ol{\si}$, $\ul{\si}$ 
are the conformally invariant
tractor extensions of $\si$ given by \nn{ulol}.
\end{thm}

\begin{proof}
First let us compute the actions
$\om_c\sharp\sharp$ and $\ka_c\sharp\sharp$ on $\ol{\si}$ and
$\ul{\si}$. The result is
\begin{align*}
  \om_c\sharp\sharp \ol{\si}_{\vec{a}^{k-l} \vec{B}^l} 
  &= -\frac{1}{2} l(k-l)(n-k+1) 
     \X_{B^1 \dot{\vec{B}}^l}^{\quad \dot{\vec{b}}^l} 
     C_{c\ a^1}^{\ p\;\ q} \si_{p\dot{\vec{a}}^{k-l} q\dot{\vec{b}}^l} \\
  \ka_c\sharp\sharp \ol{\si}_{\vec{a}^{k-l} \vec{B}^l}
  &= -\frac{1}{2} l(l-1)(n-k+1)
     \X_{B^1 \dot{\vec{B}}^l}^{\quad \dot{\vec{b}}^l}
     C_{c\ b^2}^{\ p\;\ q} \si_{\vec{a}^{k-l} qp\ddot{\vec{b}}^l} \\
  \om_c\sharp\sharp \ul{\si}_{\vec{a}^{k+l} \vec{B}^l}
  &= \frac{1}{2} l(k+1)
     \X_{B^1 \dot{\vec{B}}^l}^{\quad \dot{\vec{b}}^l} \Bigl[
     (l-1) C_{ca^{k+2}b^2a^{k+1}} \bg_{\ddot{\vec{a}}^{k,l} \ddot{\vec{b}}^l} 
     \si_{\vec{a}^k} \\
  &  \qquad\qquad\qquad\qquad\ 
     + k C_{ca^{k+1}a^1}^{\qquad\;\ p} \bg_{\dot{\vec{a}}^{k,l} 
\dot{\vec{b}}^l}
     \si_{p\dot{\vec{a}}^k} \Bigr] \\
  \ka_c\sharp\sharp \ul{\si}_{\vec{a}^{k+l} \vec{B}^l}
  &= - \frac{1}{2} l(l-1)(k+1)
     \X_{B^1 \dot{\vec{B}}^l}^{\quad \dot{\vec{b}}^l}
     C_{ca^{k+2}b^2a^{k+1}} \bg_{\ddot{\vec{a}}^{k,l} \ddot{\vec{b}}^l}
     \si_{\vec{a}^k} {}_.
\end{align*}
Now the value $y=\frac{2}{k}$ follows immediately from Proposition \ref{charac}.
In the case of $\ol{\si}$, we can reformulate Proposition \ref{charac}
in the following way: $\si$ is a solution of (\ref{ke}) if and only if
\begin{align*}
  \na_c \ol{\si}_{\vec{a}^{k-l}\vec{B}^l}
  \stackrel{{\{c\vec{a}^{k-l}\}_0}}{=}
  \frac{l(n-k+1)}{n-k} \X_{B^1\dot{\vec{B}}^l}^{\quad \dot{\vec{b}}^l} \Bigl[
  & (k-l) C_{c\ a^1}^{\ p\ q} \si_{p\dot{\vec{a}}^{k-l}q\dot{\vec{b}}^l} \\
  & + (l-1) C_{c\ b^2}^{\ p\ q} \si_{\vec{a}^{k-l}qp\dot{\vec{b}}^l} 
    \Bigr],
\end{align*}
cf. (\ref{skew_long}). Thus the value $x=\frac{2}{n-k}$ follows.
\end{proof}

\begin{rem} Note that the connections \nn{dcoup} preserve the SO$(p,q)$ 
  symmetry type (over tensor indices) and SO$(p+1,q+1)$ symmetry type
   of the any tensor-tractor field they
  act on. The coupled tractor-Levi Civita connection $\nd$ has this property.
  Then the $\om_c\sharp\sharp$ action preserves these symmetries since
  $\om_c$ is a 1-form taking values in the tensor product of orthogonal
  tractor endomorphisms tensor with orthogonal tensor endomorphisms. Similarly 
$\kappa_c$ is a 1-form taking values in the tensor square of orthogonal
  tractor endomorphisms.

Note also that the action $C_{ab}\sharp$ of the Weyl tensor on tensors may
in a natural way be viewed as a conformal action of the tractor
curvature $\Om_{ab}\sharp$ on {\em tensors}. (For example contract
each tensor index ``$c$'' into a $Z_C{}^c$ and then apply the usual
action of $\Om_{ab}\sharp$ on these tractor indices. Finally remove
each the new tractor index by contracting with $Z^C{}_{e}$. The result is
conformally invariant since $\Om_{ab}{}^C{}_DX^D=0$.) If we extend the
action $\Omega_{ab}\sharp$ to tensors in this way, then the connections
$\nd^x$ and $\nd^y$ become simply ${\na}_c^x = \na_c + x
\ka_{c}\sharp\sharp$ and ${\na}_c^y = \na_c + y\ka_{c}\sharp\sharp$
with $x$ and $y$ as above.
\end{rem}

\section{Applications: Helicity raising and lowering and almost 
Einstein manifolds} \label{almosthell}

In the first part here we will assume the structure is almost Einstein
in the sense of \cite{Go-Prague}. This is a manifold with a conformal
structure and a section $\al \in \mathcal{E}[1]$ satisfying $\left[
\na_{(a}\na_{b)_0} + P_{(ab)_0} \right] \al =0$. Equivalently there is
a standard tractor $I_A$ that is parallel with respect to the normal
tractor connection $\na$. It follows that $I_A := \frac{1}{n}D_A \al =
Y_A \al +Z_A^a \na_a \al -\frac{1}{n} X_A (\De+P) \al$, 
for some section $\al\in \ce[1]$,
and so
$X^AI_A=\al$ is non-vanishing on an open dense subset of $M$ and on
this subset $g=\al^{-2}\bg$ is an Einstein metric (where, recall $\bg$
is the conformal metric). In particular any conformally Einstein
manifold is almost Einstein but in general the converse is not true.

In this setting we immediately have the Theorem which follows.  Recall
that in a particular choice of metric a $k$-form $\si$ is a Killing
form if it is a solution of \nn{ifv} with $\tau$ identically 0. Let us
term a $k$-form $\si$ a {\em dual-Killing form} if it is a solution of
\nn{ifv} where instead $\rho$ is identically 0 (since on oriented
manifolds the Hodge dual of a Killing form is a dual-Killing form and
vice versa).
\begin{thm} \label{t.einstein}
Let us consider a $k$--form $\si_{\form{a}^k} \in \ce^k[k+1]$. Then, for 
$k\in \{1,\cdots ,n\}$,
$$
\begin{aligned}
  \ol{\ol{\si}}_{\vec{a}^{k-1}} 
  :&= \al \na^p \si_{\vec{a}^{k-1}p} 
      - (n-k+1) (\na^p \al) \si_{\vec{a}^{k-1}p} &\in& \, \ce^{k-1}[k] 
\end{aligned}
$$
is conformally invariant. For $ k\in \{0,\cdots ,n-1\}$,
$$
\begin{aligned}
  \ul{\ul{\si}}_{\vec{a}^{k+1}} 
  :&= \al \na_{a^{k+1}} \si_{\vec{a}^k} 
      - (k+1) (\na_{a^{k+1}} \al)  \si_{\vec{a}^k} &\in& \, \ce^{k+1}[k+2]
\end{aligned}
$$
is conformally invariant. If $\si$ is a solution of (\ref{ke})
then we have the following equivalences:
%$\ol{\ol{\si}}$ and $ \ul{\ul{\si}}$ satisfy the following:
\begin{equation} \label{curvcs}
\begin{split}
\na_{\{c} \ol{\ol{\si}}_{\vec{a}^{k-1}\}_0} &=0  \ekv 
    C_{ca^1}^{\quad\ pq} \si_{\dot{\vec{a}}^{k-1}pq} \stackrel{{\{c\vec{a}^{k-1}\}_0}}{=} 0 \\
  \na_{\{c}^{} \ul{\ul{\si}}_{\vec{a}^{k+1}\}_0} &=0  \ekv
    C_{ca^{k+1}a^1}^{\qquad\quad p} \si_{\dot{{\vec{a}}}^k p} 
\stackrel{{\{c\vec{a}^{k+1}\}_0}}{=} 0
\end{split}
\end{equation}
for $2\leq k\leq n-1$ and $1\leq k \leq n-2$, respectively.  In the case
that the first curvature condition is satisfied then the
corresponding conformal Killing form $\ol{\ol{\si}}_{\vec{a}^{k-1}}$
is a Killing form away from the zero set of $\alpha$, and in the
Einstein scale $g=\alpha^{-2}\bg$. In the case that the second curvature 
condition is satisfied then the
corresponding conformal Killing form $\ul{\ul{\si}}_{\vec{a}^{k-1}}$
is a dual-Killing form away from the zero set of $\alpha$, and in the
Einstein scale $g=\alpha^{-2}\bg$.
\end{thm}
\begin{proof}
  The first part of the proposition follows from relations
  $\ol{\ol{\si}}_{\vec{a}^{k-1}} = I^B \ol{\si}_{\vec{a}^{k-1}B}$ and
  $\ul{\ul{\si}}_{\vec{a}^{k+1}} = I^B \ul{\si}_{\vec{a}^{k+1}B}$
  where the forms $\ol{\si}_{\vec{a}^{k-1}B}$ and
  $\ul{\si}_{\vec{a}^{k+1}B}$ are defined by \nn{ulol} in Section
  \ref{s.helicity}.  The result \nn{curvcs} follows from Proposition
  \ref{charac} and continuity, since the tractor $I^B$ is parallel and
  $I^BX_B$ is non-vanishing on an open dense set in the manifold.  For
  the final points note that, from the formulae for
  $\ol{\ol{\si}}_{\vec{a}^{k-1}}$ and $\ul{\ul{\si}}_{\vec{a}^{k+1}}$ given
  in the first part of the theorem, it is clear that these are,
  respectively, coclosed and closed in the Einstein scale
  $g=\alpha^{-1}\bg$ given off the zero set of $\alpha$.
\end{proof}

\noindent {\bf Remarks:} 1.\ Note that the first curvature condition on the
right-hand side of \nn{curvcs} is that $(C \blacklozenge \si)=0$.
That is that the projection of $C\sharp \si$ to $\ce(1,k-1)[k-1]$
should vanish everywhere. Similarly the second is simply that the
(unique up to a multiple) projection of $C\sharp \si$ to
$\ce(1,k+1)_0[k+1]$ should vanish everywhere. Note that in the case
that the manifold is oriented then the second curvature condition is
exactly that the Hodge dual of $\si$ satisfies the first condition (as
applied to $(n-k)$-form solutions of \nn{ke}).

2.\ Note that on an almost Einstein manifold with a conformal Killing
$k$-form such that $(C \blacklozenge \si)=0$ then, according to the
Theorem, on the open dense set where $\alpha$ is
non-vanishing there is a scale so that $\ol{\ol{\si}}$ is a Killing
form. But the section $\alpha$ does not necessarily give a global
metric whereas the form $\ol{\ol{\si}}$ is a globally defined conformal
Killing form. A similar comment applies to $\ul{\ul{\si}}$.
  
\begin{cor} If $\si_{ab}$ is a conformal Killing 2-form then
  $$
  \ol{\ol{\si}}_{a} = \al \na^p \si_{ap} - (n-1) (\na^p \al)
  \si_{ap} $$
  is a  conformal Killing vector field (i.e.\ solution of
  (\ref{ke}) with $k=1$). 
If  $\si'_{\form{a}^{n-2}}$ is a conformal Killing $(n-2)$-form then
$$
\begin{aligned}
  \ul{\ul{\si}}'_{\vec{a}^{n-1}} 
  :&= \al \na_{a^{n-1}} \si'_{\form{a}^{n-2}} 
      - (n-1) (\na_{a^{n-1}} \al)  \si'_{\form{a}^{n-2}} &\in& \, \ce^{n-1}[n]
\end{aligned}
$$
is a conformal Killing $(n-1)$-form. 
 Away from the zero set of 
$\alpha$, $ \ol{\ol{\si}}_{a}$ is a Killing vector for the Einstein metric 
$g=\alpha^{-2}\bg$, while in this scale $\ul{\ul{\si}}'_{\vec{a}^{n-1}}$ is a dual-Killing form.
\end{cor}
\begin{proof}
This is just the Theorem above for $k=2$. The condition
$C_{(ab)_0}^{\quad\ \ pq} \si_{pq}$ is trivially satisfied, and, hence, so 
too is the dual condition (cf.\ point 1. of the Remark above). 
\end{proof} 
\noindent Note that a weaker form of the first part of the Corollary has 
been proved 
(by a direct computation)
 in \cite[7.2]{Sem-hab}.\\
\noindent {\bf Remark:}
 Note that according to the Corollary, on Einstein 4-manifolds a
 non-parallel conformal Killing 2-form implies the existence of either
 a non-trivial Killing vector field or a non-trivial dual-Killing 3-form.
 Thus if the 4-manifold is also oriented then, in any case,  a
 non-parallel conformal Killing 2-form determines a non-trivial
 Killing vector field.

\vspace{1ex}

The first part of the theorem is valid also for $k=1$ in
the sense, that if $\si_a$ satisfies (\ref{ke}) then $\ol{\ol{\si}} :=
\al \na^p \si_p - n (\na^p \al) \si_p$ is (conformally invariant and)
another almost Einstein scale.
This is easily seen as follows. Let us write
$\si_{CD}:={\Bbb D}_{CD}^a\si_a$, where $\D$ was defined for Lemma
\ref{D_invariance}. Then
\begin{equation} \label{grad_si_CD}
  \na_a \si_{CD} = \Om^p_{\ aCD}{} \si_p~.
\end{equation}
by Lemma \ref{l.nonclosed_noninv}. 
Note that $I^D\si_{CD}$ is parallel with respect
to the normal tractor connection $\na$ since
$$
\na_a I^D{\Bbb D}_{CD}^a\si_a = (\na_a \si_{CD}) I^D = \si^p
\Om_{paCD} I^D =0 . $$
Then the result follows from Theorem 3.1
of \cite{GoNur} since $\ol{\ol\si}= X^CI^D \si_{CD}$.

Some related results follow. Following \cite{GoNur} we term a metric (or
conformal structure) {\em weakly generic} if the Weyl curvature is
injective as bundle map $TM\to\otimes^3 TM$.
\begin{prop} 
  (i) If $\si_a$ is a non-homothetic conformal Killing vector field
  (i.e.\ a $k=1$ solution of (\ref{ke}) with non-constant $\nd_a
  \si^a$) on an Einstein manifold then there exists a non-trivial
  conformal gradient field. That is a non-trivial  solution $\tilde{\si}_a$ of
  (\ref{ke}) which is exact for the Einstein scale.\\
(ii) If a weakly generic
 conformally Einstein manifold $M$ admits a conformal Killing vector
 field $\si^a$, then $\si^a$ is a homothety for any Einstein metric in
 the conformal class. 
\end{prop}

\begin{proof}
  Let us write $I^1_D:=I_D$ and $I^2_C:=\si_{CP} I^P$, where
$\si_{CP}= {\Bbb D}_{CP}^a\si_a$.  These parallel tractors determine a
parallel tractor 2-form tractor $I^1_{[C}I^2_{D]}$. Let us write
$\tilde{\si}_a:= \frac{1}{2} \X_{\;\ a}^{CD} I^1_{[C}I^2_{D]}$. (Note
that from the last part of Theorem \ref{main} it follows immediately
that $\tilde{\si}_a $ is a conformal Killing field hence 
$\Om_{\ aCD}^p\tilde{\si}_p=0$ by (\ref{grad_si_CD}). 
Thus
$C_{abc}^{\quad p} \tilde{\si}_p^{}=0$.)
  
Since $I^1_D$ and $ I^2_C$ are parallel and the top slot of
$I^2_C$ is $\ol{\ol\si}= X^CI^D \si_{CD}$ it follows (Theorem 3.1 of
\cite{GoNur}) that $I^2_C=\frac{1}{n}D_C\ol{\ol{\si}}$. To compute 
$\tilde{\si}_a$ let us write explicitly
\begin{align*}
  I_D^1 &= Y_D \al + Z_D^d \na_d \al - \frac{1}{n} X_D(\De + J) \al \\
  I_C^2 &= Y_C \ol{\ol{\si}}+ Z_C^c \na_c \ol{\ol{\si}} 
           - \frac{1}{n} X_C(\De + J) \ol{\ol{\si}}.
\end{align*}
Here we have used the formula \nn{Dform} for the tractor $D$ operator. Now 
it follows easily that
$\tilde{\si}_a$ is $(\na_a \al) \ol{\ol{\si}} - \al (\na_a \ol{\ol{\si}})$
up to a (nonzero) scalar multiple. 
(From this formula, it is also easy to verify by
a direct computation that $\tilde{\si}_a$ satisfies (\ref{ke}).) In 
the Einstein scale $\alpha$ we have $\nd \alpha=0$, whence 
$\tilde{\si}_a= -\nd_a(\al \ol{\ol{\si}}) = 
-\nd_a(\al^2 \na^p \si_p)$.\\

(ii) This is an immediate consequence of part (i) since a conformal it
is well known (and an easy exercise to verify) that any conformal
gradient field $\tilde{\si}^a$ necessarily satisfies $C_{ab}{}^c{}_p 
\tilde{\si}^p=0$.
\end{proof}

One can easily access further results along these lines, but manifolds
admitting a conformal gradient field are rather well understood and
there are many classification results due to, for example, H.W.\
Brinkman, J.P.\ Bourguignon, D.V.\ Alekseevskii and others. For some
recent progress and indication of the state of art see \cite{KRII}.

Theorem \ref{t.einstein} exploited the standard tractor $I_A$ which
(corresponds to an almost Einstein scale $\al$ and) is parallel with
respect to the normal tractor connection $\na$.  Here we drop the
assumption that the manifold is almost Einstein and assume instead
that the manifold is equipped with a conformal Killing field
$\si^a$. Then we use the tractor $\si_{AB} := \mathbb{D}_{AB}^p \si_p$ 
(given by (\ref{si_A}))
provided by the conformal Killing form $\si_a$. This is not, in
general, parallel with respect to the normal tractor connection
$\na$. Rather, we obtained (\ref{grad_si_CD}) in Lemma 
\ref{l.nonclosed_noninv}.

\begin{thm} \label{t.ckf}
For each pair $\si \in \ce^1[2]$ and $\ta \in \ce^k[k+1]$ 
$$
\begin{aligned}
  \check{\ol{\ta}}_{\vec{a}^{k-2}} 
  :&= 2 \si^p \na^q \ta_{\vec{a}^{k-2}pq}
      + (n-k+1) (\na^p \si^q) \ta_{\vec{a}^{k-2}pq} \quad k\in\{2,\cdots,n\}
      \end{aligned}
$$
is a conformally invariant section of $\ce^{k-2}[k-1]$, and 
$$
\begin{aligned}
  \check{\ul{\ta}}_{\vec{a}^{k+2}}
  :&= 2 \si_{a^{k+1}} \na_{a^{k+2}} \ta_{\vec{a}^k} 
      + (k+1) (\na_{a^{k+1}} \si_{a^{k+2}}) \ta_{\vec{a}^k}, 
  \quad k\in\{0,\cdots,n\!-\!2\}
\end{aligned}
$$
is a conformally invariant section of $\ce^{k+2}[k+3]$.
If $\si$ and $\ta$ are 
solutions of (\ref{ke}) then the following is satisfied: for $3\leq k\leq n-1$ 
$\check{\ol{\ta}}_{\vec{a}^{k-2}}$, is a solution of (\ref{ke}) 
if and only if
$$ (n-k+1) C_{\ c}^{r\ pq} \ta_{\vec{a}^{k-2}pq} \si_r^{}
   + (k-2) C_{ca^1}^{\quad pq} \ta_{p\dot{\vec{a}}^{k-2}qr} \si^r 
   \stackrel{\{c\vec{a}^{k-2}\}_0}{=} 0 $$
and, for $1\leq k\leq n-3$, $\check{\ul{\ta}}_{\vec{a}^{k+2}}$, is a solution 
of (\ref{ke}) 
if and only if
$$ 2 C_{ca^{k+1}a^1}^{\qquad\quad p} \ta_{p \dot{{\vec{a}}}^k} \si_{a^{k+2}}
   - C_{\ ca^{k+1}a^{k+2}}^{p} \ta_{\vec{a}^k} \si_p^{}
   \stackrel{\{c\vec{a}^{k+2}\}_0}{=} 0. $$
\end{thm}

\begin{proof}
The first part of the proposition follows from relations
$\check{\ol{\ta}}_{\vec{a}^{k-2}} = \ol{\ta}_{\vec{a}^{k-2}RS} \si^{RS}$ and
$\check{\ul{\ta}}_{\vec{a}^{k+2}} = \ul{\ta}_{\vec{a}^{k+2}RS} \si^{RS}$.
The second part is a result of a direct computation. Using 
Proposition \ref{charac} and (\ref{grad_si_CD}) we obtain the following:
\begin{align*}
\na_c & \ol{\ta}_{\vec{a}^{k-2}RS} \si^{RS} \stackrel{\{c\vec{a}^{k-2}\}_0}{=}
    (\na_c \ol{\ta}_{\vec{a}^{k-2}RS}) \si^{RS} +
    \ol{\ta}_{\vec{a}^{k-2}RS} \na_c \si^{RS} \\
  & \stackrel{\{c\vec{a}^{k-2}\}_0}{=}
    \frac{2(n-k+1)}{n-k} \X_{RS}^{\ \ s} 
    \left[ (k-2) C_{c\ a^1}^{\ p\;\ q} \ta_{p\dot{\vec{a}}^{k-2}qs}
           -  C_{c\ s}^{\ p\ q} \ta_{p\dot{\vec{a}}^{k-2}qa^1} \right] \si^{RS} \\
  & \qquad\quad + \ta_{\vec{a}^{k-2}RS} \Om^{p\ RS}_{\ c} \si_p \\
  & \stackrel{\{c\vec{a}^{k-2}\}_0}{=}
    \frac{n\!-\!k\!+\!1}{n\!-\!k} \left[ 
    (n\!-\!k\!+\!1) C_{\ c}^{s\ pq} \ta_{\vec{a}^{k-2}pq} \si_s^{}
 + (k\!-\!2) C_{ca^1}^{\quad pq} \si_{p\dot{\vec{a}}^{k-2}qs} \si^s \right], \\
\na_c & \ul{\ta}_{\vec{a}^{k+2}RS} \si^{RS} \stackrel{\{c\vec{a}^{k+2}\}_0}{=}
    (\na_c \ul{\ta}_{\vec{a}^{k+2}RS}) \si^{RS} +
    \ul{\ta}_{\vec{a}^{k+2}RS} \na_c \si^{RS} \\
  & \stackrel{\{c\vec{a}^{k+2}\}_0}{=}
    - 2 (k+1) \X_{RS}^{\ \ s} C_{ca^{k+1}a^1}^{\qquad\  p} 
      \ta_{p\dot{\vec{a}}^k} \bg_{a^{k+2}s} \si^{RS}
    + \ul{\ta}_{\vec{a}^{k+2}RS} \Om^{p\ RS}_{\ c} \si_p \\
  & \stackrel{\{c\vec{a}^{k+2}\}_0}{=}
    - (k+1) \left[
    2 C_{ca^{k+1}a^1}^{\qquad\quad p} \ta_{p \dot{{\vec{a}}}^k} \si_{a^{k+2}}
    - C_{\ ca^{k+1}a^{k+2}}^{p} \ta_{\vec{a}^k} \si_p^{}  
    \right].
\end{align*}
\end{proof}
\noindent Note for the cases of a conformal Killing 3-form $\tau$ the
first curvature condition of the Theorem is satisfied by any conformal
gradient vector field $\si$.

Now it is obvious how to obtain more general results for couples of
conformal Killing forms $\si \in \cE^l[l+1]$ and $\ta \in \cE^k[k+1]$
where $1 \leq k,l \leq n-1$.  We set $\si_{\form{A}^{l+1}} :=
\mathbb{D} \si$ and define $\check{\ol{\ta}}_{\vec{a}^{k-l-1}}
:= \ol{\ta}_{\vec{a}^{k-l-1}\form{A}^{l+1}} \si^{\form{A}^{l+1}}$ and
$\check{\ul{\ta}}_{\vec{a}^{k+l+1}} :=
\ul{\ta}_{\vec{a}^{k+l+1}\form{A}^{l+1}} \si^{\form{A}^{l+1}}$ for $0
\leq k-l-1 \leq n$ and $0 \leq k+l+1 \leq n$, respectively.  The case
$l=1$ is described in the previous Theorem and in general, the obstructions
for $\check{\ol{\ta}}_{\vec{a}^{k-l-1}}$ and
$\check{\ul{\ta}}_{\vec{a}^{k+l+1}}$ to be solutions of (\ref{ke}) are
very similar to the cases $l=1$. (In the proof of these new cases, we
replace $\na_c \si^{RS}$ by $\na_c \si^{\form{A}^{l+1}}$. The latter
is, in general quite complicated but we actually need only
'$\Z$--slot' and '$\Y$--slot' which are similar to the case $l=1$.)

\begin{cor}
Let $\si_a \in \cE_a[2]$ be a solution of (\ref{ke}) and write 
$\mu_{bc} := \na_{[b}\si_{c]}$ (in a choice of scale). Then the section
$$ \si_{a^0}\mu_{a^1a^2} \cdots \mu_{a^{2p-1}a^{2p}}
   \in \cE^{2p+1}[2p+2], \ p \leq \lfloor \frac{n-2}{2} \rfloor $$
is conformally invariant. 
If $\si_{a^0}^{} C_{a^1a^2c}^{\qquad d} \si_d^{}=0$ then  
it %If the latter condition is satisfied then the displayed section
is a solution of (\ref{ke}) for any $1\leq p \leq \lfloor \frac{n-2}{2} \rfloor $.
\end{cor}

\begin{proof}
For $p=1$, this is Theorem \ref{t.ckf} applied to  $\ta := \si \in \cE^1[2]$.
If the curvature condition is satisfied then it is easily checked
that applying the same Theorem to $\si_a$ and $\ta := \si_{a^0}\mu_{a^1a^2}$, 
we obtain the case $p=2$. Repeating this procedure, the general case follows.
\end{proof} 

Let us note there are several results in \cite{Sem-MathZ} related to
those in this section, see Propositions 3.4 and 3.5 in
\cite{Sem-MathZ}.  These concern a special case satisfying that $\na_c
\mu_{a^0a^1}$ is pure trace (which implies that $\si_a$ is an
eigensection for the Schouten tensor viewed as a section of
End$(TM)$).  This immediately yields $\si_{a^0}^{} C_{a^1a^2c}^{\qquad
  p} \si_p^{}=0$ using (\ref{prolong_start}).

\vspace{1ex}

Our last application concerns conformal Killing $m$-tensors.  These
are valence $m$ symmetric trace-free tensors $t_{b\cdots c} \in
\cE_{(b\cdots c)_0}[2m]$ which are solutions of the conformally
invariant equation $\na_{(a} t_{b\cdots c)_0} =0$. Obviously, any
conformal Killing form $\si_a \in \cE_a[2]$ yields a conformal Killing
tensor $\si_{(a} \cdots \si_{b)_0}$.  Note that generalising the $m=2$
version of this observation we have the following.  If $\si_{\form{a}}
\in \cE_{\form{a}^k}[k+1]$ is conformal Killing form then
$\si_{(a}^{\;\ \dform{c}} \si_{b)_0\dform{c}}^{} \in \cE_{(ab)_0}[4]$,
is a conformal Killing 2-tensor.  (The special case of this where
$\si$ is a conformal Killing 2-form appeared in
\cite[4.1(4)]{Stepanov}.)  This follows from (\ref{prolong_start}) by
a direct computation or from the relation $\si_{(a}^{\;\ \dform{e}}
\si_{b)_0\dform{e}}^{} = \frac{1}{(n-k+1)^2} \ol{\si}_{(a}^{\;\ 
  \dform{E}} \ol{\si}_{b)_0\dform{E}}^{}$ (which holds since $X_A$ and
$Z^A_a$ are orthogonal), and Propositions \ref{charac} and
\ref{coupKill}. The point here is that one applies the normal tractor
$\na_c$ connection to $\ol{\si}_{(a}^{\;\ \dform{E}}
\ol{\si}_{b)_0\dform{E}}^{}$ to obtain $2\ol{\si}_{(a}^{\;\ \dform{E}}
\na_b^{\phantom{ \dform{E}}} \ol{\si}_{c)_0\dform{E}}^{}$ after the
projection to $\cE_{(abc)_0}[4]$.  Then from Proposition \ref{charac}
and again the orthogonality of $X_A$ and $Z^A_a$ we may replace $\na$
by $\tilde{\nd}$ to obtain $2\ol{\si}_{(a}^{\;\ \dform{E}}
\tilde{\na}_b^{\phantom{ \dform{E}}} \ol{\si}_{c)_0\dform{E}}^{}$.
But then by Proposition \ref{coupKill} the last expression vanishes.
It is clear this example generalises and so we have the following
Theorem.

\begin{thm}\label{steplike}
Suppose $\si^1,\cdots , \si^m$ is a collection of conformal Killing
forms of respective ranks $r_1,\cdots , r_m$ where $(\sum^m r_i)-m$ is
an even number. Then 
$$
\si^1_{(a}\cdot \si^2_b\cdot~\cdots~\cdot \si^m_{c)_0}
$$ is a conformal Killing $m$-tensor, where $\si^1_{a}\cdot
\si^2_b\cdot~\cdots~\cdot \si^m_c$ indicates any contraction of the collection 
$\si^1,\cdots , \si^m$
over the
suppressed indices.
\end{thm}
\noindent Of course it will often be the case that a given contraction
$\si^1_{a}\cdot \si^2_b\cdot~\cdots~\cdot \si^m_c$ vanishes upon
projection to the trace-free part. However it is easy to proliferate
non-trivial examples.

\end{document}